\documentclass[preprint, 3p]{elsarticle}

\usepackage{lipsum}
\makeatletter
\def\ps@pprintTitle{%
 \let\@oddhead\@empty
 \let\@evenhead\@empty
 \def\@oddfoot{}%
 \let\@evenfoot\@oddfoot}
\makeatother

\usepackage{amssymb}
\usepackage{amsthm}
\usepackage{amsbsy}
\usepackage{amsmath}
\usepackage{algorithm,algpseudocode}

\usepackage{xcolor}
\usepackage{makecell}
\usepackage{tabularx}

\usepackage{lineno,hyperref}
\modulolinenumbers[5]

\journal{} %{Computer Methods in Applied Mechanics and Engineering}

\begin{document}

\begin{frontmatter}

\title{POD-DL-ROM: enhancing deep learning-based reduced order models for nonlinear parametrized PDEs by proper orthogonal decomposition}

%% Group authors per affiliation:
\author[1]{Stefania Fresca}
\ead{stefania.fresca@polimi.com}
\author[1]{Andrea Manzoni\corref{cor1}}
\ead{andrea1.manzoni@polimi.com}
\cortext[cor1]{Corresponding author}

\address[1]{MOX - Dipartimento di Matematica, Politecnico di Milano, P.zza Leonardo da Vinci 32, 20133 Milano, Italy}

\begin{abstract}
Deep learning-based reduced order models (DL-ROMs) have been recently proposed to overcome common limitations shared by conventional reduced order models (ROMs)  -- built, e.g., through proper orthogonal decomposition (POD) -- when applied to nonlinear time-dependent parametrized partial differential equations (PDEs). These might be related to {\em (i)} the need to deal with projections onto high dimensional linear approximating trial manifolds, {\em (ii)} expensive hyper-reduction strategies, or {\em (iii)} the intrinsic difficulty to handle physical complexity with a linear superimposition of modes.
All these aspects are avoided when employing DL-ROMs, which learn in a non-intrusive way both the nonlinear trial manifold and the reduced dynamics, by relying on deep (e.g., feedforward, convolutional, autoencoder) neural networks. Although extremely efficient at testing time, when evaluating the PDE solution for any new testing-parameter instance, DL-ROMs require an expensive training stage, because of the extremely large number of network parameters to be estimated.
In this paper we propose a possible way to avoid an expensive training stage of DL-ROMs, by {\em (i)} performing a prior dimensionality reduction through POD, and {\em (ii)} relying on a multi-fidelity pretraining stage, where different physical models can be efficiently combined.
The proposed POD-DL-ROM is tested on several (both scalar and vector, linear and nonlinear) time-dependent parametrized PDEs (such as, e.g., linear advection-diffusion-reaction, nonlinear diffusion-reaction, nonlinear elastodynamics, and Navier-Stokes equations) to show the generality of this approach and its remarkable computational savings.

\end{abstract}

\begin{keyword}
reduced order modeling, deep learning, proper orthogonal decomposition, dimensionality reduction, parametrized PDEs
\end{keyword}

\end{frontmatter}

%\linenumbers

\section{Introduction}

Performing the numerical approximation of parametrized partial differential equations (PDEs) for multiple parameter values, or solving them in real-time, is unaffordable by means of full order, traditional high-fidelity techniques, such as the Galerkin-finite element method \cite{quarteroni1994numerical}. Replacing a full order model (FOM) by a reduced order model (ROM), featuring a much lower dimension, yet capable to express the physical features of the problem at hand, is the main goal of reduced order modeling techniques, among which the reduced basis (RB) method represents one of the most popular options \cite{quarteroni2016reduced}. The basic assumption underlying the RB method is that the solution of a parametrized PDE lies on a low-dimensional manifold, which can be approximated by a linear trial subspace spanned by a set of basis functions \cite{benner2017model, benner2015asurvey}, built from a set of FOM snapshots employing, e.g., proper orthogonal decomposition (POD). In this case, the ROM approximation is given by the linear superimposition of POD modes, whose degrees of freedom (depending both on time and parameters) result from the solution of a low-dimensional (nonlinear, dynamical) system, obtained through a (Petrov-)Galerkin projection onto a linear test subspace, which might coincide with the trial subspace. Being able to assemble such a ROM efficiently, for any new parameter instance, is granted at the price of a further {\em hyper-reduction} stage on the FOM arrays \cite{FGMQ_19}. %, usually exploiting techniques like the (discrete) empirical interpolation  \cite{barrault2004empirical,chaturantabut2010nonlinear,Maday2009383} or the gappy POD \cite{WILLCOX2006208}.

ROMs for parametrized PDEs rely on a suitable offline-online computational splitting: computationally expensive tasks required to build the low-dimensional subspace, and assemble all the ROM arrays, are performed once for all during the so-called {\em offline} (or {\em ROM training}) stage, then allowing to compute -- ideally -- in an extremely efficient way the ROM approximation for any new parameter value, during the so-called {\em online} (or {\em ROM testing}) stage.
This computational strategy, however, breaks down if {\em (i)} the dimension of the linear trial subspace becomes very large (compared to the intrinsic dimension of the solution manifold being approximated), or {\em (ii)} the hyper-reduction stage, required to approximate parameter-dependent nonlinear terms -- at the limit, the whole residual vector and, possibly, the Jacobian matrix when dealing with implicit solvers such as, e.g., the Newton method -- requires linear subspaces whose dimension becomes very large, too, in order to provide an approximation to FOM arrays sufficiently accurate.
These can be recurrent issues when dealing with nonlinear, time-dependent parametrized PDEs, and {\em (i)} we aim at building a ROM able to provide problem approximations uniformly accurate over the whole parameter space, {\em (ii)} the parametrized problem features coherent structures  (e.g., transport or wave phenomena) that propagate over time and depend on parameters. Last, but not least, ensuring ROM stability might require additional computational efforts, such as, e.g., when dealing with fluid flows using a mixed formulation (e.g., a velocity-pressure formulation for Navier-Stokes equations), see, e.g., \cite{DALSANTO2019186,doi:10.1002/nme.4772, RozzaHuynhManzoni2013}. \\

To overcome these drawbacks, we have recently proposed in \cite{fresca2020comprehensive, fresca2020deep} a strategy to construct deep learning-based ROMs (DL-ROMs) for nonlinear time-dependent parametrized PDEs in a non-intrusive way, exploiting deep neural networks as main building block, and a set of FOM snapshots. In particular, the DL-ROM  technique allows approximating both the solution manifold of the PDE  by means of a low-dimensional, nonlinear trial manifold, and the nonlinear dynamics of the generalized coordinates on such reduced trial manifold, as a function of the time coordinate and the parameters. Both {\em (i)} the nonlinear trial manifold and {\em (ii)} the reduced dynamics are learnt in a non-intrusive way, thus avoiding the projection stage typical of the RB method; in particular, the trial manifold is learnt by means of the decoder function of a convolutional autoencoder (CAE) neural network, whereas the reduced dynamics through a (deep) feedforward neural network (DFNN), and the encoder function of the CAE (see Section \ref{sec:sec3} for further details).

DL-ROMs outperform the RB method -- even involving local reduced bases -- regarding both numerical accuracy (for the same ROM dimension) and computational efficiency during the {\em online} (or {\em testing}) stage, when applied to problems that are typically challenging for the RB method (such as, e.g., linear transport equation, nonlinear diffusion-reaction equations whose solution develops moving fronts depending on parameters) or problems featuring reduced bases with usually large dimension. A key aspect, still open in the setting of DL-ROMs and which this paper is mainly devoted to, deals with the computational efficiency of DL-ROMs during the {\em offline} (or {\em training}) stage, which is also related with the curse of dimensionality. Indeed, training the two networks representing the building blocks of the DL-ROM entails a number of networks parameters to be estimated (and, correspondingly, training data dimensions) that blow up with the dimension $N_h$ of the FOM -- this latter being related with, e.g., the mesh size required by a Galerkin-finite element approximation of the PDEs. Indeed, DL-ROMs have been applied so far to the reduction of scalar PDEs, with at most $N_h = O(10^4)$ degrees of freedom.   \\ 

In this paper we propose a strategy to enhance DL-ROMs in order to make the {\em offline} training stage dramatically faster, allowing for much larger FOM dimensions, without affecting  the number of networks parameters to be estimated and, ultimately, network complexity. This strategy exploits {\em (i)} dimensionality reduction of FOM snapshots through randomized POD %, performed by means of randomized SVD
(rPOD)  \cite{szlam2004implementation}, to be considered as the action of the first {\em layer} of the CAE, rather than the way to generate the linear trial manifold, as done instead in traditional POD-Galerkin ROMs, and {\em (ii)} a suitable multi-fidelity pretraining stage \cite{goodfellow2016deep}, where different models (built, e.g., by considering coarser discretizations or simplified physical models) can be efficiently combined, to iteratively initialize the network parameters.

%Improving on the weakest aspect -- and still taking advantage  of the key properties -- of DL-ROMs, is the aim of this %work.
These substantial enhancements of the DL-ROM technique provide a new way to build deep learning-based ROMs, which we refer to POD-DL-ROM. The resulting strategy represents a suitable combination of the best features of deep learning (DL) algorithms and POD  -- namely, the non-intrusive character of the former, and the simplicity, combined with the rigorous mathematical foundation, of the latter --  at the same time fixing their weaknesses  -- namely, the curse of dimensionality of DL-ROM for increasing FOM dimensions, and the modest approximation properties of POD-based linear trial manifolds for some classes of nonlinear parametrized  PDEs. As a result, the POD-DL-ROM technique provides not only substantial gains during the {\em offline} training stage if compared to the DL-ROM strategy -- keeping fixed the dimension of the FOM that we used to generate snapshots -- but also during the {\em online} testing stage, if compared to POD-Galerkin ROMs, also in those cases where this latter technique provides satisfying results in terms of both accuracy and efficiency.  \\ % Indeed, we have assessed the computational performance of POD-DL-ROMs in {\em (i)} a linear advection-diffusion-reaction problem,  {\em (ii)} a nonlinear diffusion-reaction problem arising from cardiac electrophysiology, {\em (iii)} structural mechanics, and
%{\em (iv)} fluid dynamics. As a result, POD-DL-ROMs are shown to yield extremely efficient numerical approximations to (scalar and vector) nonlinear time-dependent parametrized PDEs, thus providing a {\em turn-key} strategy to build ROMs only relying on a set of FOM snapshots, and ultimately leading to the possibility to simulate in more than real-time, during the {\em online} testing stage, physical phenomena occurring on a time scale of seconds. \\

{The structure of the paper is as follows. In Section \ref{sec:2} we briefly recall the construction of DL-ROMs by reinterpreting the classical ideas behind linear projection-based methods for parametrized PDEs, comparing our strategy with alternative ways to build ROMs relying on machine/deep learning algorithms. In Section \ref{sec:3} we then describe the POD-DL-ROM technique showing {\em (i)} how to enhance the DL-ROM technique by means of (randomized) POD, {\em (ii)} how to rely on a multi-fidelity pretraining stage, and {\em (iii)} how to generalize this technique to the case of vector problems. In Section \ref{sec:num_res} we assess the numerical accuracy and efficiency on four different test cases, namely {\em (i)} a linear advection-diffusion-reaction problem,  {\em (ii)} a nonlinear diffusion-reaction problem arising from cardiac electrophysiology, {\em (iii)} nonlinear elastodynamics for hyperelastic compressible materials, and {\em (iv)} fluid dynamics, 
%advection-diffusion-reaction, monodomain, elastodynamics and Navier-Stokes equations, 
also showing the great versatility of the proposed technique. %Finally, we draw some conclusions.
}

\section{An overview of deep learning-based ROMs} \label{sec:2}

Before introducing the main features of POD-DL-ROMs, we review the construction of DL-ROMs and highlight the main differences between our framework and existing techniques in literature.

\subsection{DL-ROMs in a nutshell}
\label{sec:sec2}

For the sake of generality, we consider a generic nonlinear, time-dependent  PDE depending on a set of input parameters $\boldsymbol{\mu} \in \mathcal{P}$, where  the parameter space $\mathcal{P} \subset \mathbb{R}^{n_{\boldsymbol{\mu}}}$ is given by a bounded and closed set; input parameters might represent physical or geometrical properties of the system, like, e.g., material properties, initial and boundary conditions, or the shape of the domain.  Even if we only consider physical parameters, handling geometrical parameters does not require additional efforts.
We adopt a fully algebraic perspective and assume to start from the high-fidelity (spatial) approximation of the PDE, which we refer to as full order model (FOM). Regardless of the spatial discretization adopted -- such as, e.g., the finite element method, Isogeometric Analysis  or the spectral element method -- the FOM can be expressed as a nonlinear parametrized dynamical system. Hence, given $\boldsymbol{\mu} \in \mathcal{P}$, we aim at solving the initial value problem \vspace{-0.1cm}
\begin{equation}
\label{eq:FOM}
\begin{cases}
{\bf M}(\boldsymbol{\mu}) \mathbf{\dot{u}}_h(t;\boldsymbol{\mu}) = \mathbf{f}(t, \mathbf{u}_h(t;\boldsymbol{\mu}); \boldsymbol{\mu}) \quad t \in (0, T),\\
\mathbf{u}_h(0;\boldsymbol{\mu})=\mathbf{u}_0(\boldsymbol{\mu}), \vspace{-0.1cm}
\end{cases}
\end{equation}
where:
\begin{itemize}
\item[-] $\mathbf{u}_h:[0,T) \times \mathcal{P} \rightarrow \mathbb{R}^{N_h}$ is the parametrized  solution of (\ref{eq:FOM});
\item[-] $\mathbf{u}_0 : \mathcal{P} \rightarrow \mathbb{R}^{N_h}$ is the initial datum;
\item[-] $\mathbf{f} : (0,T) \times  \mathbb{R}^{N_h}  \times \mathcal{P} \rightarrow \mathbb{R}^{N_h}$ is a (nonlinear) function, encoding the system dynamics;
\item[-] ${\bf M}(\boldsymbol{\mu}) \in \mathbb{R}^{N_h \times N_h}$ is the parametric, mass matrix of this parametric FOM;  without any loss of generality, ${\bf M}(\boldsymbol{\mu})$ is assumed here to be a symmetric positive definite matrix.
\end{itemize}
The FOM dimension $N_h$  is related with the finite dimensional subspaces introduced for the sake of space discretization -- here $h>0$ denotes a discretization parameter, such as the maximum diameter of elements in a computational mesh; consequently, $N_h$ can be extremely large if the PDE problem describes complex physical behaviors and/or high degrees of accuracy are required to its solution. In order to solve problem (\ref{eq:FOM}), suitable time discretizations are employed, such as backward differentiation formulas (BDFs) \cite{quarteroni1994numerical} and generalized-$\alpha$ \cite{chung1993atime} methods. We thus aim at approximating, in an efficient way, the set
\begin{equation}
\mathcal{S}_h = \{ \mathbf{u}_h(t ; \boldsymbol{\mu} ) \; | \; t \in [0, T) \; , \; \boldsymbol{\mu} \in \mathcal{P} \subset \mathbb{R}^{n_{\mu}} \} \subset \mathbb{R}^{N_h},
\label{eq:solution_manifold}
\end{equation}
of solutions to problem (\ref{eq:FOM}) when $(t ; \boldsymbol{\mu} )$ varies in $[0, T) \times \mathcal{P} $, also referred to as {\em solution manifold}.
If for any $\boldsymbol{\mu} \in \mathcal{P}$, problem (\ref{eq:FOM}) admits a unique solution,  for each $t \in (0,T)$, the intrinsic dimension of the solution manifold is at most $n_{\boldsymbol{\mu}} + 1  \ll N_h$, where $n_{\boldsymbol \mu}$ is the number of parameters (time plays the role of an additional coordinate). In this case, each element $\mathbf{u}_h(t ; \boldsymbol{\mu} ) \in \mathcal{S}_h$ can be described in terms of at most $n_{\boldsymbol{\mu}} + 1$ intrinsic coordinates, even if their explicit characterization is, in practice, computationally unaffordable. Equivalently, the tangent space to the manifold at any given $\mathbf{u}_h(t ; \boldsymbol{\mu} )$ is spanned by $n_{\boldsymbol{\mu}} + 1$ basis vectors. \\

When dealing with traditional projection-based ROMs, the approximated solution to  \eqref{eq:FOM} is sought in the subspace  {$\textnormal{Col}(\mathbf{V}_n)$} of dimension $n \ll N_h$, spanned by the $n$ columns of \textcolor{black}{$\mathbf{V}_n \in \mathbb{R}^{N_h \times n}$}. Hence, a linear ROM looks for an approximation $\mathbf{\tilde{u}}_h(t; \boldsymbol \mu) \approx \mathbf{u}_h(t; \boldsymbol \mu)$ of the form 
\begin{equation}
\label{linear_reconstructed_solution}
\mathbf{\tilde{u}}_h(t; \boldsymbol \mu) = {\bf V}_n\mathbf{u}_n(t; \boldsymbol{\mu}) , %+ \mathbf{u}_0(\boldsymbol{\mu})
\end{equation}
where %$\mathbf{u}_0(\boldsymbol{\mu}) \in \mathbb{R}^{N_h}$ and
$\mathbf{\tilde{u}}_h : [0,T) \times \mathcal{P} \rightarrow \textcolor{black}{\tilde{\mathcal{S}}_h^{n,lin}}$. Hence,  the \textit{reduced linear  trial manifold} is given by
\begin{equation}
\label{manifold_linear}
\textcolor{black}{\tilde{\mathcal{S}}_h^{n,lin} = \{
{\bf V}_n\mathbf{u}_n(t; \boldsymbol{\mu})  \ | \mathbf{u}_n(t; \boldsymbol{\mu}) \in \mathbb{R}^n,    \ t \in [0, T)  \; , \; \boldsymbol{\mu} \in \mathcal{P} \subset {\mathbb{R}}^{n_{\mu}} \} \subset \mathbb{R}^{N_h}.}
\end{equation}

In  POD-Galerkin ROMs,  \textcolor{black}{$\tilde{\mathcal{S}}_h^{n,lin}$} is spanned by the first $n$ singular vectors of 
\begin{equation}\label{eq:snap_matrix1}
{\bf S} = [   \mathbf{u}_h^1(t^1; \boldsymbol \mu_1) \; | \; \ldots \; | \; \mathbf{u}_h^1(t^{N_t}; \boldsymbol \mu_1) \; | \; \ldots \; | \; \ldots
\mathbf{u}_h^1(t^1 ; \boldsymbol \mu_{N_{train}}) \; | \; \ldots \; | \; \mathbf{u}_h^1(t^{N_t} ; \boldsymbol \mu_{N_{train}}) ],
\end{equation}
a matrix collecting FOM solutions (or {\em snapshots}) computed for different parameter values  $\boldsymbol{\mu}_1, \ldots, \boldsymbol{\mu}_{N_{train}} \in \mathcal{P}$, suitably sampled\footnote{Sampling frequency in time does not necessarily coincide with the time stepping rule used to discretize in time  (\ref{eq:FOM}); similarly, singular value decomposition can be done either on ${\bf S}$, or following a two-stage procedure, operating SVD across samples in time for each $\boldsymbol \mu_i$, $i=1, \ldots, N_{train}$, and then on the resulting collection of selected singular vectors.   } over the parameter space, at different time instants $\{t^1, \ldots, t^{N_t}\} \subset [0,T]$. 

In the POD-Galerkin case, to model the reduced dynamics of the system, we replace $\mathbf{u}_h(t; \boldsymbol{\mu})$ by (\ref{linear_reconstructed_solution}) in system (\ref{eq:FOM}), and impose that the residual
\begin{equation}
\label{residual}
\mathbf{r}_h({\bf V}_n\mathbf{u}_n (t; \boldsymbol{\mu})) = {\bf M}(\boldsymbol{\mu}) {\bf V}_n \mathbf{\dot{u}}_n(t; \boldsymbol{\mu}) - \mathbf{f}(t, {\bf V}_n\mathbf{u}_n (t; \boldsymbol{\mu}); \boldsymbol{\mu})
\end{equation}
is  orthogonal to \textcolor{black}{$\tilde{\mathcal{S}}_h^{n,lin}$}. This condition yields the following POD-Galerkin ROM
\begin{equation}
\label{ROM}
\begin{cases}
{\bf M}_n(\boldsymbol{\mu}) \mathbf{\dot{u}}_n(t;\boldsymbol{\mu}) = \mathbf{f}_n(t, \mathbf{V}_n\mathbf{u}_n(t;\boldsymbol{\mu}); \boldsymbol{\mu}) \quad t \in (0, T),\\
\mathbf{u}_n(0;\boldsymbol{\mu})= {\bf V}_n^T \mathbf{u}_0(\boldsymbol{\mu}),
\end{cases}
\end{equation}
where:
\begin{itemize}
\item[-] ${\bf M}_n(\boldsymbol{\mu}) = {\bf V}_n^T {\bf M}(\boldsymbol{\mu}) {\bf V}_n$ is the {\em reduced} mass matrix;
\item[-] $\mathbf{f}_n(t, \mathbf{V}_n\mathbf{u}_n(t;\boldsymbol{\mu}); \boldsymbol{\mu}) = {\bf V}_n^T \mathbf{f}(t, {\bf V}_n\mathbf{u}_n(t;\boldsymbol{\mu}); \boldsymbol{\mu})$;
\item[-] $\mathbf{u}_n(0;\boldsymbol{\mu}) = {\bf V}_n^T \mathbf{u}_h(0;\boldsymbol{\mu})$ is the initial condition for $\mathbf{u}_n(t;\boldsymbol{\mu})$ associated with the initial condition for $\mathbf{u}_h(t;\boldsymbol{\mu})$, where we have assumed,
 without loss of generality, that ${\bf V}_n^T {\bf V}_n = {\bf I} \in \mathbb{R}^{n \times n}$.
\end{itemize}
The two main bottlenecks often arising with POD-Galerkin ROMs are {\em (i)} the increasing dimension $n \gg n_{\mu} + 1$ of  the low-dimensional POD subspaces, much larger than the intrinsic dimension of the solution manifold, and {\em (ii)} the need to rely on hyper-reduction techniques  to assemble the operators appearing in the ROM (\ref{ROM}) in order not to  rely on expensive $N_h$-dimensional arrays  \cite{FGMQ_19}.

DL-ROMs have been introduced in \cite{fresca2020comprehensive} to overcome these  limitations of POD-Galerkin ROMs, and further applied to cardiac electrophysiology in \cite{fresca2020deep}. A DL-ROM describes both the  trial manifold and the reduced dynamics (corresponding to the matrix ${\bf V}_n$ and the projection stage, respectively, in the POD-Galerkin case) through deep  neural networks, which are trained on a set of FOM snapshots. In this way, DL-ROMs completely avoid the {\em projection} stage, are non-intrusive, and can be cheaply evaluated once trained.  In particular:
 \begin{itemize}
%We denote by $N_{train}$ and $N_{test}$ the number of training and testing parameter instances, respectively; the ROM dimension is again denoted by $n  \ll N_h$.
%
\item   to describe the system dynamics on a suitable reduced nonlinear trial manifold (a task which we refer to as {\em reduced dynamics learning}), the intrinsic coordinates of the ROM approximation are defined as
\begin{equation}
\label{eq:reduced_dynamics}
{\mathbf{u}}_n(t; \boldsymbol \mu, \boldsymbol{\theta}_{DF}) = {\boldsymbol{\phi}}_n^{DF}(t; \boldsymbol \mu, \boldsymbol{\theta}_{DF}), % \vspace{0.15cm}
\end{equation}
where ${\boldsymbol{\phi}}_n^{DF}(\cdot; \cdot, \boldsymbol{\theta}_{DF}) : {\mathbb{R}}^{(n_{\mu} +1)} \rightarrow {\mathbb{R}}^n$ is a DFNN, consisting in the subsequent composition of a nonlinear activation function,  applied to a linear transformation of the input, multiple times \cite{goodfellow2016deep}.  Here $\boldsymbol{\theta}_{DF}$ denotes the vector of {parameters} of the DFNN, collecting all the corresponding weights and biases of each layer of the DFNN;

\item  to model the reduced nonlinear trial manifold  $\tilde{\mathcal{S}}_h^n \approx \mathcal{S}_h$ (a task which we refer to as  {\em reduced trial manifold learning}) we employ the decoder function of a CAE  \cite{lecun1998gradient, hinton1994autoencoders}, that is, % More precisely, $\tilde{\mathcal{S}}_h$ takes the form
\begin{equation}
\label{eq:nonlinear_manifold}
%\begin{aligned}
\textcolor{black}{\tilde{\mathcal{S}}_h^n =  \{ {\mathbf{f}}^D_h(\mathbf{u}_n(t; \boldsymbol{\mu}, {\boldsymbol{\theta}_{DF}}); \boldsymbol{\theta}_{D})
\; | \; \mathbf{u}_n(t; \boldsymbol{\mu}, \boldsymbol{\theta}_{DF}) \in {\mathbb{R}}^{n},    \ t \in [0, T)  \; \textnormal{and} \; \boldsymbol{\mu} \in \mathcal{P} \subset {\mathbb{R}}^{n_{\mu}} \} \subset \mathbb{R}^{N_h},}
%= \\
%& \{ \mathbf{f}_{h}^D(\boldsymbol{\phi}_n^{DF}(t; \boldsymbol{\mu}, {\boldsymbol{\theta}_{DF}}); \boldsymbol{\theta}_{D})
%\; | \; \boldsymbol{\phi}_n^{DF}(t; \boldsymbol{\mu}, \boldsymbol{\theta}_{DF}) \in \mathbb{R}^{n}, \\
%& \ t \in [0, T)  \; \textnormal{and} \; \boldsymbol{\mu} \in \mathcal{P} \subset \mathbb{R}^{n_{\mu}} \} \subset \mathbb{R}^{N},
%\end{aligned}
\end{equation}
where ${\mathbf{f}}^D(\cdot; {\boldsymbol{\theta}}_D) : {\mathbb{R}}^n \rightarrow {\mathbb{R}^{N_h}}$ denotes the   decoder function of a CAE obtained as   the composition of several layers (some of which are convolutional), depending upon a vector ${\boldsymbol{\theta}}_{D}$ collecting all the corresponding weights and biases. % of each layer of the decoder.
\end{itemize}

The DL-ROM approximation $\mathbf{\tilde{u}}_h(t; \boldsymbol \mu) \approx \mathbf{u}_h(t; \boldsymbol \mu)$  is then given by
\begin{equation}
\mathbf{\tilde{u}}_h(t; \boldsymbol \mu, \theta_{DF}, \theta_D) = {\mathbf{f}}_h^D({\boldsymbol{\phi}}_n^{DF}(t; \boldsymbol{\mu}, {{\boldsymbol{\theta}}_{DF}}); \boldsymbol{\theta}_{D}).
\label{eq:reconstructed_solution}
\end{equation}
The encoder function ${\mathbf{f}}_{n}^E( \cdot ; \boldsymbol{\theta}_{E})$ of the convolutional AE, provided when carrying out its training on the FOM snapshots, can then be exploited to map the FOM solution associated to $(t, \boldsymbol{\mu})$ onto a low-dimensional representation
\begin{equation}
\label{eq:encoder}
 {\mathbf{\tilde{u}}_n}(t; \boldsymbol{\mu}, \boldsymbol{\theta}_{E}) = {\mathbf{f}}_{n}^E(\mathbf{u}_h(t; \boldsymbol{\mu}); \boldsymbol{\theta}_{E});
\end{equation}
 $\mathbf{f}_{n}^E(\cdot; \boldsymbol{\theta}_E) : \mathbb{R}^{N_h} \rightarrow \mathbb{R}^n$ denotes  the encoder function, depending upon a vector ${\boldsymbol{\theta}}_E$ of parameters.

Computing the DL-ROM approximation of ${\mathbf{u}}_h(t; \boldsymbol{\mu}_{test})$, for any possible $t \in (0,T)$ and $\boldsymbol{\mu}_{test} \in \mathcal{P}$, corresponds to the testing  stage of a DFNN and of the decoder function of a convolutional AE; this does not require the evaluation of the encoder function. The training stage consists in solving the following optimization problem, in the  variable $\boldsymbol{\theta} = (\boldsymbol{\theta}_{E}, \boldsymbol{\theta}_{DF}, \boldsymbol{\theta}_{D})$, after the snapshot matrix $\mathbf{S}$ has been formed:
\begin{equation}
\min_{\boldsymbol{\theta}} \mathcal{J}(\boldsymbol{\theta}) = \min_{\boldsymbol{\theta}} \frac{1}{N_s}\sum_{i=1}^{N_{train}} \sum_{k=1}^{N_t} \mathcal{L}(t^k, \boldsymbol \mu_i; \boldsymbol{\theta}),
\label{eq:minimization_problem}
\end{equation}
where {$N_s = N_{train} N_t$} and
\begin{equation}
\mathcal{L}(t^k, \boldsymbol{\mu}_i;  {\boldsymbol{\theta}}) = \frac{\omega_h}{2} \mathcal{L}_{rec}(t^k, \boldsymbol{\mu}_i;  {\boldsymbol{\theta}}) + \frac{1-\omega_h}{2}       \mathcal{L}_{int}(t^k, \boldsymbol{\mu}_i;  {\boldsymbol{\theta}}),
\label{eq:loss_encoder}
\end{equation}
where
\[
\mathcal{L}_{rec}(t^k, \boldsymbol{\mu}_i;  {\boldsymbol{\theta}}) = \| \mathbf{u}(t^k; \boldsymbol{\mu}_i) - \mathbf{\tilde{u}}(t^k; \boldsymbol{\mu}_i,  {\boldsymbol{\theta}_{DF}, \boldsymbol{\theta}_D})\|^2,
\]
\[
\mathcal{L}_{int}(t^k, \boldsymbol{\mu}_i;  {\boldsymbol{\theta}}) =  \| \tilde{\mathbf{u}}_n(t^k; \boldsymbol{\mu}_i, \boldsymbol{\theta}_E) -  {\mathbf{u}}_n(t^k; \boldsymbol{\mu}_i,  {\boldsymbol{\theta}_{DF}})\|^2
\]
with $\omega_h \in [0,1]$. The \emph{per-example} loss function (\ref{eq:loss_encoder}) combines the reconstruction error (that is, the error between the FOM solution and the DL-ROM approximation) and the error between the  {intrinsic coordinates} and the output of the encoder.

\subsection{The content of this paper and other existing approaches}

Several recent works have shown possible applications of DL algorithms for solving PDEs -- thanks to their ability of effectively approximating nonlinear maps, and by their ability to learn from data and generalize to unseen data --  both from a theoretical \cite{kutyniok2019atheoretical, opschoor2020deep, laakmann2020efficient} %, petersen2018optimal}
and a computational standpoint. Regarding this latter, we mention, for instance, physics-informed neural networks (PINNs)
%For example, Karniadakis and coauthors replace the FOMs by physics-informed neural networks (PINNs)
\cite{raissi2018hidden,raissi2017physics1,raissi2017physics2,raissi2019physics} or %trained by minimizing the residual of the PDEs, computed by means of automatic differentiation \cite{baydin2018automatic}. In particular, in \cite{raissi2018deep} this framework is applied to the Kuramoto-Sivashinsky equation, in a chaotic regime, however, the algorithm leads to unsatisfactory approximations. A similar approach to PINNs can be found in \cite{han2017solving} and \cite{yang2018physics} where
physics-informed deep generative models \cite{han2017solving,yang2018physics}. % that is generative models based on DL algorithms, aiming at computing the solution of PDEs, are presented. DL techniques for parametrized PDEs have also been proposed in various recent works.
DL algorithms and artificial neural networks (ANN), such as feedforward neural networks, are becoming more and more popular in reduced order modeling, too. In particular:
\begin{itemize}
\item  ANNs have been employed to model the reduced dynamics in a data-driven and less intrusive way (avoiding, e.g., the costs entailed by the projection stage of projection-based ROMs); for instance, in  \cite{guo2018reduced,guo2019data, hestaven2018non-intrusive, kast2020anon, wang2019non} the use of ANNs or Gaussian Processes (GPs) regression models has been proposed to approximate the mapping from the input parameters (and, possibly, time) to the reduced coefficients, as an alternative to the assembly and solution of the reduced order system arising from POD-Galerkin ROMs, however still using a linear trial manifold built through POD; similar strategies have also been introduced in \cite{kani2017dr-rnn,mohan2018adeep,wan2018data, pulch2020machine, berzins2020standardized}. A hybrid strategy proposed in \cite{chena2020physics} merges ANN-based regression models and PINNs,  training the network by minimizing the mean squared residual error of the reduced order equation on a set of points in parameter space; similar results can be found in \cite{kani2018reduced}. Moreover, feedforward and recurrent neural networks have been exploited in   \cite{san2018neural, wang2020recurrent} to address closure problems and model the effects of discarded modes on the set of retained POD modes.  {Very recently, an ANN-based methodology is proposed to learn mappings between infinite-dimensional spaces for parametric PDEs  in \cite{bhattacharya2020model}, exploiting POD and showing mesh-independent properties.}

\item ANNs have been used to describe the reduced trial manifold where the approximation is sought (thus avoiding the linear superimposition of POD modes), either relying on a minimum residual formulation to derive the ROM \cite{carlberg2018model, kim2020afast}, or without considering an explicit parameter dependence in the differential problem that is considered \cite{gonzalez2018deep}. For instance, projection-based ROMs are built in   \cite{carlberg2018model} by performing a projection of the FOM onto a nonlinear trial manifold identified by means of the decoder function of a CAE -- hence, still requiring the assembling and the solution of a ROM as in traditional POD-Galerkin ROMs. The use of CAEs has also been proposed in \cite{gonzalez2018deep}, where a reduced  trial manifold is generated through  a deep convolutional recurrent AE, which is then used to train a Long Short-Term Memory (LSTM) neural network that models the reduced dynamics.  

\end{itemize}

%\textcolor{red}{In this paper we propose a DL framework aiming at constructing non-intrusive ROMs to address nonlinear time-dependent parametrized PDEs combining the best features of some of the approaches mentioned above. In particular, we employ deep feedforward neural networks (DFNNs) and convolutional autoencoders (CAEs) to  construct the reduced trial manifold and model the reduced dynamics on it. In particular, {\em (i)} a DFNN is used to map $(t,\boldsymbol{\mu})$  }

%%%%%%%

Our DL-ROM approach \cite{fresca2020comprehensive, fresca2020deep}  combines and  improves the techniques introduced in \cite{gonzalez2018deep,carlberg2018model}. As detailed in Section \ref{sec:sec2}, the {\em (i)} nonlinear trial manifold is learnt by using the decoder function of a CAE, while {\em (ii)} the dynamics on the reduced manifold is modeled through a DFNN and the encoder function of a CAE. These two tasks are achieved simultaneously, by training both the CAE and the DFNN network architectures at the same time, by minimizing a loss function weighting two terms -- see \eqref{eq:loss_encoder} -- one for each task. The resulting procedure thus  %In this respect, we are able to design a flexible framework capable to handle parameters affecting both PDE operators and data, which 
 avoids both the expensive projection stage of \cite{carlberg2018model} and the training of a more expensive LSTM network \cite{gonzalez2018deep}. 
%In the DL-ROM, the intrusive construction of a ROM is replaced by the evaluation of the ROM generalized coordinates through a DFNN taking  only $(t, \boldsymbol{\mu})$ as inputs. 
Moreover, the DL-ROM technique is purely data-driven, non-intrusive: it only relies on the computation of a set of FOM snapshots. In this respect, DL does not replace the high-fidelity FOM as, e.g., in the works by Karniadakis and coauthors \cite{raissi2018hidden,raissi2017physics1,raissi2017physics2,raissi2019physics,raissi2018deep}; rather, DL techniques are built upon it, to enhance the repeated evaluation of the FOM for different values of the parameters.  The computational benefits %\\
%
%  ultimately yielding ROMs whose dimension is nearly equal (if not equal) to  the intrinsic dimension of the solution manifold that we aim at approximating.  
%
%
%\textcolor{red}{The proposed DL-ROM can efficiently provide solutions to time-dependent nonlinear parametrized problems, especially if compared to common linear (projection-based) ROMs, built for instance through a POD-Galerkin RB method.
%In particular, the benefits 
introduced by the use of DL-ROMs can be summarized as follows:
\begin{itemize}
\item the dimension of the DL-ROM can be kept extremely small, very close (or even equal) to the dimension of the solution manifold $n_{\mu} + 1$;
\item the DL-ROM can be queried at any desired time instant, without requiring the solution of a dynamical system until that time, differently from projection-based ROMs such as, e.g, POD-Galerkin ROMs;
\item the time resolution required by the DL-ROM can be chosen to be larger than the one required by the numerical solution of dynamical systems at hand (see, e.g., \cite{fresca2020deep});
\item DL-ROMs avoid the use of (intrusive and very often extremely expensive) hyper-reduction techniques, which are instead required by POD-Galerkin ROMs; 
\item DL-ROMs can avoid to account for those auxiliary variables of a problem which we might not be interested  into (as pressure, compared to velocity, in fluid flow problems, or the gating variables, compared to the electric potential, in cardiac electrophysiology % dynamics of varia\textit{(i)} to handle efficiently those terms that depend nonlinearly on either the solution or the input parameters, and \textit{(ii)} to account for the dynamics of variables, solution of the problem at hand, we are not interested in reconstructing, as, for instance, in the case of cardiac electrophysiology problems DL-ROMs allow to approximate the electric potential without the need of computing the gating variables 
\cite{fresca2020deep}).
\end{itemize}
{For all these reasons, DL-ROMs %completely avoid both the assembly and the projection stages typical of traditional projection-based ROMs, thus 
tremendously improve the computational efficiency of ROMs during the {\em online} testing stage. 
%In addition, outputs of clinical interest, such as ACs and APs, can be more efficiently evaluated by the DL-ROM technique than by a FOM built, for example, through the  FE method, while maintaining a high level of accuracy.
However, the {\em offline} training stage of DL-ROM would still depend on $N_h$, a fact which ultimately might entail overwhelming training times and costs when $N_h$ is moderately large.} 

The POD-DL-ROM technique proposed in this work, thanks to a prior dimensionality reduction relying on POD and a suitable multi-fidelity pretraining, greatly enhances the efficiency of the DL-ROM during the training phase, thus dramatically decreasing training computational times, as shown by the numerical results discussed in following Sections.

\section{A new deep learning-based reduced order model}\label{sec:3}

POD-DL-ROMs provide a new, general-purpose, ROM approach combining a data dimensionality reduction obtained through POD % (equivalent to principal component analysis in statistics \cite{hastie2001theelements}, or Karhunen-Lo\`eve expansion in stochastic applications), 
%a widespread technique used to compute a reduced trial manifold in linear ROMs, 
with the DL-ROM approach %, a technique exploiting DL models to learn the reduced nonlinear trial manifold and the reduced dynamics on it, thus enabling the construction of ROMs characterized by efficient testing computational times
 \citep{fresca2020comprehensive, fresca2020deep}. After introducing the POD-DL-ROM approach, we discuss in more details some of its building blocks, the extension to vector problems, finally reporting detailed algorithms for the   {\em offline} (or training) and the {\em online} query (or testing) stages.

\subsection{POD-enhanced DL-ROMs (POD-DL-ROMs)}

The POD-DL-ROM technique consists in applying the DL-ROM technique to the intrinsic coordinates of a linear trial manifold generated through randomized singular value decomposition (rSVD) and approximating $\mathcal{S}_h$; alternatively, it can be seen as a ROM technique in which a two-step  dimensionality reduction is performed: first, POD (realized through rSVD) is applied on a set of FOM snapshots, then a DL-ROM  \textcolor{black}{is built to approximate the map between $(t, \boldsymbol{\mu})$ and the POD generalized coordinates. In this way, all the DL-ROM features allowing its very efficient testing time are retained. As a matter of fact,} as shown in Section \ref{sec:sec2}, DL-ROMs might imply overwhelming training costs (and times) when the FOM dimension $N_h$ becomes moderately large,  although remaining extremely efficient at testing time. We emphasize that very often ROMs are designed to be efficient only regarding their {\em online} performances, no matter how expensive is the offline stage. Our (more ambitious) goal is instead to realize a ROM able of efficient computational performance during both offline  and online stages, compared to classical projection-based ROMs.  \\ 

Using randomized POD (see Section \ref{sec:sec3}), we first  build the  $N$-dimensional \textcolor{black}{subspace} $\textnormal{Col}(\mathbf{V}_N)$ spanned by the $N \leq N_h$ columns of ${\bf V}_N \in \mathbb{R}^{N_h \times N}$, the matrix of the first $N$ singular vectors of the snapshot matrix ${\bf S}$. We denote  the dimension of the linear manifold  by $N$, to distinguish it from the dimension $n$ of the nonlinear trial manifold, and to emphasize that this dimension can be taken (much) larger with respect to the one of the reduced linear trial manifold that would have been exploited in a POD-Galerkin ROM. Indeed, here linear dimensionality reduction is performed only for the sake of data compression, to avoid to feed training data of dimension $N_h$.

\textcolor{black}{The POD-DL-ROM approximation $\mathbf{\tilde{u}}_h(t; \boldsymbol{\mu}, \boldsymbol{\theta}_{DF}, \boldsymbol{\theta}_{D})$ of the FOM solution ${\mathbf{u}}_h(t; \boldsymbol{\mu})$ is given by
\begin{equation*}
\mathbf{\tilde{u}}_h(t; \boldsymbol{\mu}, \boldsymbol{\theta}_{DF}, \boldsymbol{\theta}_{D}) = \mathbf{V}_N \mathbf{\tilde{u}}_N(t; \boldsymbol{\mu}, \boldsymbol{\theta}_{DF}, \boldsymbol{\theta}_{D}), 
\end{equation*}
that is, it is sought in a linear trial manifold\footnote{Equivalently, we have replaced the original solution manifold $\mathcal{S}_h$ with the $N$-dimensional linear manifold
\begin{equation*}
\label{eq:manifold_S_N}
\textcolor{black}{
\mathcal{S}_{N}^{h} = \{ \mathbf{V}_N^T\mathbf{u}_{h}(t ; \boldsymbol{\mu} ) \; | \; t \in [0, T) \; \textnormal{and} \; \boldsymbol{\mu} \in \mathcal{P} \subset \mathbb{R}^{n_{\mu}} \} \subset \mathbb{R}^{N}}.
\end{equation*}
As in the case of $\mathcal{S}_h$, the  intrinsic dimension of $\mathcal{S}_{N}^{h}$ is at most $n_{\mu} + 1 \ll N$.} of (potentially large) dimension $N$,
%
% 
% 
%
% Indeed, the POD-DL-ROM approximation of  is sought in the linear trial manifold
\begin{equation}
\label{manifold_linear_N}
\tilde{\mathcal{S}}_h^{N,lin} = \{
{\bf V}_N \mathbf{\tilde{u}}_N(t; \boldsymbol{\mu}, \boldsymbol{\theta}_{DF}, \boldsymbol{\theta}_{D})  \ | \mathbf{\tilde{u}}_N(t; \boldsymbol{\mu}, \boldsymbol{\theta}_{DF}, \boldsymbol{\theta}_{D}) \in \mathbb{R}^N,    \ t \in [0, T)  \; \textnormal{and} \; \boldsymbol{\mu} \in \mathcal{P} \subset {\mathbb{R}}^{n_{\mu}} \} \subset \mathbb{R}^{N_h},
\end{equation}
by applying the DL-ROM strategy  of Section \ref{sec:sec2} to approximate $\mathbf{V}_N^T\mathbf{u}_{h}(t ; \boldsymbol{\mu} )$ -- rather than ${\mathbf{u}}_h(t; \boldsymbol{\mu})$. 
The DL-ROM approximation 
$\mathbf{\tilde{u}}_N(t; \boldsymbol{\mu}, \boldsymbol{\theta}_{DF}, \boldsymbol{\theta}_{D}) \approx \mathbf{V}_N^T\mathbf{u}_{h}(t ; \boldsymbol{\mu} )$ takes the form
\begin{equation}
\mathbf{\tilde{u}}_N(t; \boldsymbol \mu, \boldsymbol{\theta}_{DF}, \boldsymbol{\theta}_D ) = {\mathbf{f}}^D_N({\boldsymbol{\phi}}_n^{DF}(t; \boldsymbol{\mu}, {{\boldsymbol{\theta}}_{DF}}); \boldsymbol{\theta}_{D}),
\label{eq:u_N_approx}
\end{equation}
and is sought in a reduced nonlinear trial manifold $\tilde{\mathcal{S}}_N^n$ of dimension $n \ll N$. By adapting the DL-ROM formulation of Section \ref{sec:sec2} to the case at hand, we have that: 
\begin{itemize}
\item to describe the system dynamics on % we aim at learning the dynamics of the set ${\mathcal{S}}_N$ on 
 the nonlinear trial manifold $\tilde{\mathcal{S}}_N^n$ -- with   $n$  as close as possible to $n_{\boldsymbol{\mu}} + 1$ --  the intrinsic coordinates of the approximation $\mathbf{\tilde{u}}_N$ are defined as  %in terms of minimal coordinates 
%by means of a DFNN; indeed, we 
\begin{equation*}
\label{eq:phi_n}
\mathbf{u}_n(t; \boldsymbol \mu) = \boldsymbol{\phi}_n^{DF}(t; \boldsymbol \mu, \boldsymbol{\theta}_{DF}),
\end{equation*}
where   $\boldsymbol{\phi}_n(\cdot ;  \cdot, \boldsymbol{\theta}_{DF}) : [0, T) \times \mathbb{R}^{n_{\mu} +1} \rightarrow \mathbb{R}^n$ is a DFNN; %, 
%\textcolor{blue}{where $\boldsymbol{\theta}_{DF}$ denotes  the vector of {parameters} of the DFNN}; %, collecting all the corresponding weights and biases of each layer of the DFNN;
\item to model the  reduced nonlinear trial manifold \textcolor{black}{$\tilde{\mathcal{S}}_N^n$}, we  employ the decoder function of a CAE, that is, 
%we define the function in (\ref{eq:Psi_N}) as
%\begin{equation*}
%\label{eq:f_N}
%\boldsymbol{\Psi}_h(\mathbf{u}_n(t; \boldsymbol{\mu}, {\boldsymbol{\theta}_{DF}}); \boldsymbol{\theta}_{D}) = \mathbf{f}_{N}^D(\mathbf{u}_n(t; \boldsymbol{\mu}, {\boldsymbol{\theta}_{DF}}); \boldsymbol{\theta}_{D}).
%\end{equation*}
%\textcolor{blue}{where $\boldsymbol{\theta}_D$ denotes the vector of parameters of the convolutional and dense layers of the decoder function.}
\begin{equation}
\begin{split}
\tilde{\mathcal{S}}_N^n = \{ \mathbf{\tilde{u}}_{N}(t; \boldsymbol \mu) = & {\mathbf{f}}^D_N({\boldsymbol{\phi}}_n^{DF}(t; \boldsymbol{\mu}, {{\boldsymbol{\theta}}_{DF}}); \boldsymbol{\theta}_{D}) \; | \\
& \; \mathbf{u}_n(t; \boldsymbol{\mu},  {{\boldsymbol{\theta}}_{DF}}) \in \mathbb{R}^{n}, \ t \in [0, T) \; , \; \boldsymbol{\mu} \in \mathcal{P} \subset \mathbb{R}^{n_{\mu}} \} 	\subset \mathbb{R}^N,
\end{split}
\label{eq:manifold_tilde_S_n}
\end{equation}
where ${\mathbf{f}}^D_N( \cdot ; \boldsymbol{\theta}_{D}) : 	\mathbb{R}^n \rightarrow \mathbb{R}^N$. %is the encoder function of a CAE.
\end{itemize}}

The encoder function  of the convolutional AE can then be exploited to map the intrinsic coordinates $\mathbf{V}_N^T \mathbf{u}_h (t, \boldsymbol{\mu})$ associated to $(t, \boldsymbol{\mu})$ onto a low-dimensional representation
\begin{equation*}
\label{eq:f_n}
 {\mathbf{\tilde{u}}_n}(t; \boldsymbol{\mu}, \boldsymbol{\theta}_{E}) = {\mathbf{f}}_{n}^E(\mathbf{V}_N^T \mathbf{u}_h(t; \boldsymbol{\mu}); \boldsymbol{\theta}_{E}),
\end{equation*}
where $\mathbf{f}_{n}^E$ denotes  the encoder function, depending upon a vector ${\boldsymbol{\theta}}_E$ of parameters.

\begin{figure}[b!]
\centering
\includegraphics[scale=0.245]{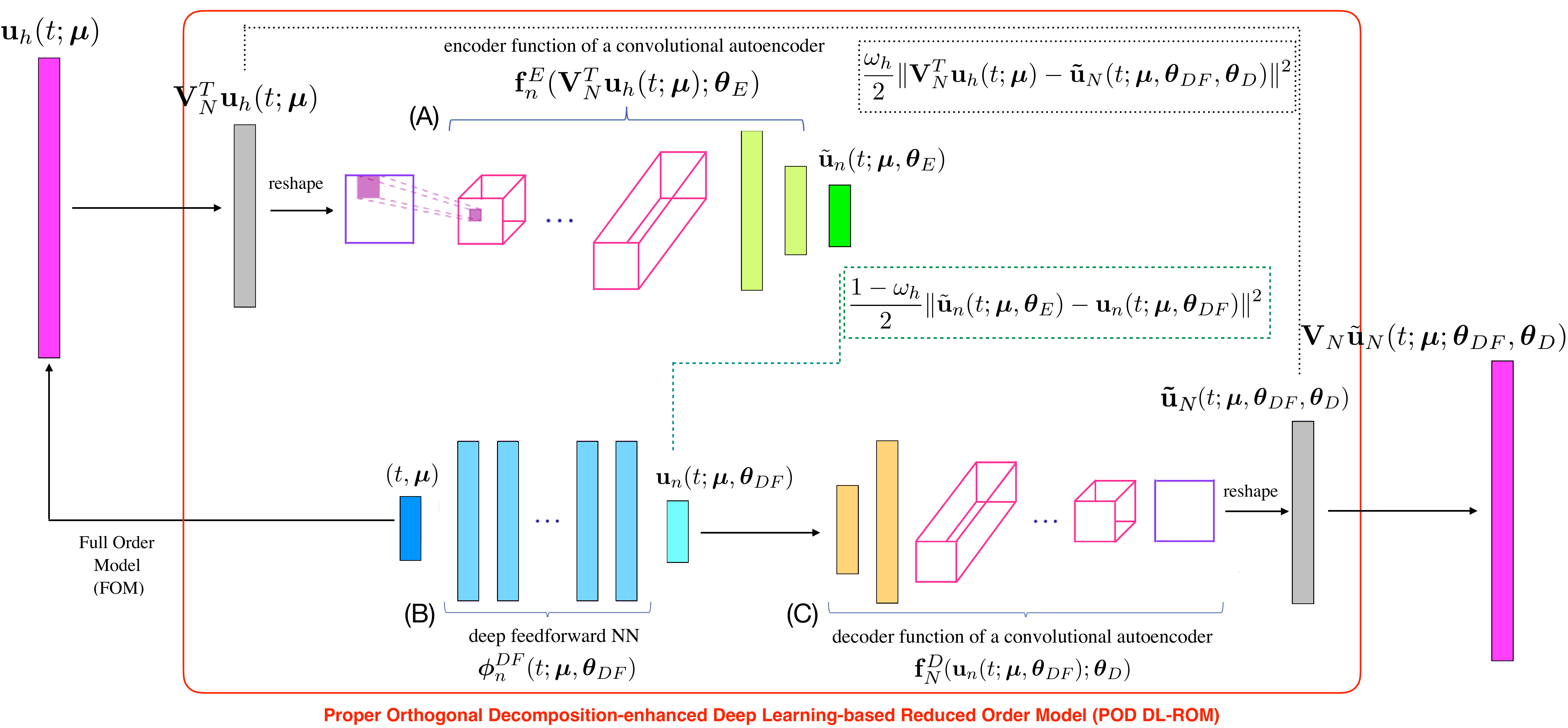}
\caption[POD-DL-ROM architecture.]{POD-DL-ROM architecture. Starting from the FOM solution $\mathbf{u}_h(t; \boldsymbol{\mu})$, the intrinsic coordinates $\mathbf{V}_N^T \mathbf{u}_h(t; \boldsymbol{\mu})$ are computed, by means of rSVD, and the neural network provides as output $\mathbf{\tilde{u}}_N(t; \boldsymbol{\mu})$, an approximation of them. The reconstructed solution $\mathbf{\tilde{u}}_h(t; \boldsymbol{\mu})$ is then recovered trough the rPOD basis matrix.}
\label{fig:Fig5-0}
\end{figure}

The architecture of the POD-DL-ROM neural network, employed at training time, is the one shown in \figurename~\ref{fig:Fig5-0}; note that,  at testing time, as in the DL-ROM technique we can discard the encoder function.

 Computing the ROM approximation (\ref{eq:u_N_approx}) by means of  a POD-DL-ROM thus consists in solving the optimization problem (\ref{eq:minimization_problem}) where  the per-example loss function \eqref{eq:loss_encoder} is now replaced by
\begin{equation}
\label{eq:loss_N}
\begin{split}
\mathcal{L}(t^k, \boldsymbol{\mu}_i;  {\boldsymbol{\theta}})  = \frac{\omega_h}{2}\| \mathbf{V}_N^T\mathbf{u}_h(t^k; \boldsymbol{\mu}_i) - \mathbf{\tilde{u}}_N(t^k; \boldsymbol{\mu}_i,  {\boldsymbol{\theta}_{DF}, \boldsymbol{\theta}_D})\|^2 \\
 + \frac{1-\omega_h}{2}  \| \tilde{\mathbf{u}}_n(t^k; \boldsymbol{\mu}_i, \boldsymbol{\theta}_E) -  {\mathbf{u}}_n(t^k; \boldsymbol{\mu}_i,  {\boldsymbol{\theta}_{DF}})\|^2,
\end{split}
\end{equation}
where
\[
\mathcal{L}_{rec}(t^k, \boldsymbol{\mu}_i;  {\boldsymbol{\theta}}) = \| \mathbf{V}_N^T\mathbf{u}_h(t^k; \boldsymbol{\mu}_i) - \mathbf{\tilde{u}}_N(t^k; \boldsymbol{\mu}_i,  {\boldsymbol{\theta}_{DF}, \boldsymbol{\theta}_D})\|^2.
\]

We point out that the POD-DL-ROM technique relies on a further level of dimensionality reduction compared to \cite{hestaven2018non-intrusive, wang2019non}, since the dimension of the reduced linear problem -- that is, the number of intrinsic coordinates included in $\mathbf{\tilde{u}}_N$ -- is decreased until (almost, or exactly) matching the intrinsic dimension $n_{\mu} + 1$ of the parametrized problem. In addition, we introduce the use of convolutional layers, which results better suited to high-dimensional spatial data, thus implying lower computational costs, with respect to dense layers. 
Moreover, the POD-DL-ROM approach allows to model and approximate the entire intrinsic coordinates vector $\mathbf{\tilde{u}}_N$ all at once, without requiring additional SVDs, if compared to the data-driven RB method, employing GPs as regression models,  proposed in  \cite{guo2019data} and further extended in  \cite{kast2020anon}. 

We highlight that by shaping the POD-DL-ROM neural network as a zero-extended neural network and assuming that the input-output map is locally Lipschitz, it is possible to prove the convergence of the framework here proposed by following the approach presented in \cite{bhattacharya2020model}.

%a data-driven RB method has been proposed in  \cite{guo2019data}, employing GPs to build a regression model for the non-intrusive approximation of the intrinsic coordinates associated to a set of POD modes. In particular, regression functions have been defined as  linear combinations of  tensor products of two GPs,  depending on either time or input  parameters; a further SVD approximation is then used to decompose the training data associated to each element of the intrinsic coordinates vector into several time- and parameter-modes, so that, ultimately, GPs are trained to approximate the modes of each intrinsic coefficient. This approach has been extended in \cite{kast2020anon}, where a multi-fidelity setup is used both for constructing a reduced space and for   GP regression. However, when dealing with high POD dimensions, the previous strategies may entail the training of a very high number of regression models. The POD-DL-ROM approach, instead, allows us to model and approximate the entire intrinsic coordinates vector all at once, without requiring additional SVDs.

\subsection{(Randomized) POD for dimensionality reduction}
\label{sec:sec3}

%POD exploits the SVD of the snapshot matrix $\mathbf{S}$ to compute an orthonormal basis, which is optimal in a least-squares sense \cite{quarteroni2016reduced}. Due to the linear superimposition of modes assumption, POD is not able to compute a reduced subspace of dimension equal (or close) to the dimension of the solution manifold, that is to say $n_{\mu} + 1$, while preserving an acceptable degree of accuracy with respect to the FOM solution, when dealing with problems featuring coherent structures propagating over time \cite{fresca2020comprehensive,fresca2020deep}. %we showed, by computing the optimal-POD solution, that is the projection of the FOM solution onto the linear trial manifold generated through POD, or the POD-Galerkin ROM approximation, that the dimension $N$ of the linear trial manifold built by means of POD, must be much higher than $n_{\mu} + 1$ in order to ensure the same level of accuracy obtained with a DL-ROM. 
Even though POD is not able to generate linear subspaces whose dimension is close to (or matches) the intrinsic dimension of the problems under consideration, it is still able to perform a moderate dimensionality reduction, thus yielding a linear subspace of dimension $N \ll N_h$. %This latter feature is exploited in the POD-DL-ROM strategy.
However, computing the SVD of the snapshot matrix $\mathbf{S} \in \mathbb{R}^{N_h \times N_{train}N_t}$ can be extremely time consuming for large-scale problems;  whenever dealing with FOMs of moderately large dimension $N_h$ and/or a high number of training-parameter or time instances, the computational time and  memory required by SVD may become prohibitive,  scaling superlinearly in %. For every large values of $N_h$ and $N_{train}N_t$, the time and the memory required by the SVD are superlinear in 
$N_h$ and $N_{train}N_t$ \cite{drineas2006fast}. 
In order to speed up computations, we rely on the randomized matrix approximation techniques developed in \cite{halko2011finding}; in particular, we exploit rSVD, which computes an approximated SVD, using randomization.

More precisely, a Gaussian random matrix $\mathbf{\Omega} \in \mathbb{R}^{N_{train}N_t \times m}$ is drawn, where $N \le m \le N_{train}N_t$; $m - N$ is called the oversampling parameter. The matrix
\begin{equation*}
\mathbf{Y} = (\mathbf{S}\mathbf{S}^T)^q\mathbf{S} \mathbf{\Omega}
\end{equation*}
is then assembled, by normally setting $q=1$ or 2; then,  an iterative QR factorization  $\mathbf{Y} = \mathbf{QR}$, is computed, where $\mathbf{Q} \in \mathbb{R}^{N_h \times m}$. Computing $\mathbf{Q}$ in such a way consists of applying an adaptive randomized range finder algorithm to approximate the range of $\mathbf{S}$ by means of a matrix $\mathbf{Q}$ whose columns are orthonormal, i.e. $\mathbf{S} \approx \mathbf{Q} \mathbf{Q}^T \mathbf{S}$ \cite{halko2011finding}. Once the matrix $\mathbf{Q}$ has been computed, it is restricted to the first $N$ columns, denoted by $\mathbf{Q} = \mathbf{Q} ( \, : \, , 1: N )  \in \mathbb{R}^{N_h \times N}$, and the SVD of the matrix
\begin{equation}
\label{rSVD_B}
\mathbf{B} =  \mathbf{Q}^T\mathbf{S} = \mathbf{\tilde{V} \tilde{\Sigma} \tilde{Z}}
\end{equation}
is computed. The SVD factorization of $\mathbf{S}$ is then recovered by setting
\begin{equation}
\label{rSVD_VN}
\mathbf{V}_N = \mathbf{Q \tilde{V}}.
\end{equation} 
The rSVD approximation can then be computed through the following steps:
\begin{enumerate}
\item an approximated basis for Col($\mathbf{S}$), i.e. $\mathbf{Q} \in \mathbb{R}^{N_h \times m}$, is obtained by using randomization,
  \item the SVD of the matrix $\mathbf{B} \in \mathbb{R}^{N \times N_{train}N_t}$ in (\ref{rSVD_B}) is computed,
  \item the matrix $\mathbf{V}_N \in \mathbb{R}^{N_h \times N}$ is recovered by means of $\mathbf{Q}$ as in (\ref{rSVD_VN}).
  \end{enumerate}
We refer to $\mathbf{V}_N$ as to the rPOD basis matrix. For further details about rSVD we refer to \cite{halko2011finding, szlam2004implementation}.

%By applying rSVD, we build a reduced linear trial manifold $\tilde{\mathcal{S}}_{h} = \textnormal{Col}(\mathbf{V})$ of dimension $N \ll N_h$, approximating the solution manifold $\mathcal{S}_h$. The ROM approximation is then provided by
%\begin{equation*}
%{\mathbf{u}}_h(t; \boldsymbol{\mu}) \approx \mathbf{\tilde{u}}_h(t; \boldsymbol{\mu}) = \mathbf{V} \mathbf{u}_{N}(t; \boldsymbol{\mu}),
%\end{equation*}
%where $\mathbf{\tilde{u}}_h : [0,T) \times \mathcal{P} \rightarrow \tilde{\mathcal{S}}_{h}$.
%Here $\mathbf{u}_{N}(t; \boldsymbol{\mu}) \in \mathbb{R}^{N}$, for each $t \in [0,T)$ and $\boldsymbol{\mu} \in \mathcal{P}$, denotes an approximation of the vector of the intrinsic coordinates of the ROM approximation on the linear subspace generated through rSVD, that is
%\begin{equation*}
%\mathbf{u}_N(t; \boldsymbol{\mu}) \approx \mathbf{V}^T\mathbf{u}_h(t; \boldsymbol{\mu}).
%\end{equation*}
%We finally define the whole set
%\begin{equation}
%\label{eq:manifold_S_N}
%\mathcal{S}_{N} = \{ \mathbf{V}^T\mathbf{u}_{h}(t ; \boldsymbol{\mu} ) \; | \; t \in [0, T) \; \textnormal{and} \; \boldsymbol{\mu} \in \mathcal{P} \subset \mathbb{R}^{n_{\mu}} \} \subset \mathbb{R}^{N},
%\end{equation}
%of the intrinsic coordinates, when $(t, \boldsymbol{\mu})$ varies in $ [0,T) \times \mathcal{P}$,  which is a, possibly nonlinear, manifold of dimension $n_{\mu} + 1 \le$ dim($\mathcal{S}_N$) $\ll N$.

\subsection{Pretraining}\label{sec:pretraining}

Directly training a model to solve a specific task can be very demanding if the model is complex and hard to optimize, and/or the task is very difficult. More viable options are either {\em (i)} to train a simpler model to solve the task, then make the model more complex, or {\em (ii)} to train the model to solve a simpler task, then move on to the final task. Both these strategies are known in the DL literature as \textit{pretraining} \cite{goodfellow2016deep}. In particular, a finer tuning of a pretrained model is equivalent to \textit{transfer learning}, if the data used to perform fine tuning are of different nature with respect to the data used during pretraining. Pretraining can then be seen as a form of transfer learning, where a pretrained model is used as initial state of the network \cite{taylor2009transfer}; this strategy works extremely well in many objects classification tasks \cite{yosinki2014how} and natural language processing problems \cite{devlin2018bert}. In the area of scientific ML, pretraining has been used, for instance, in \cite{haghighat2020adeep} where a pretrained neural network has been used to perform \textcolor{black}{parameter} identification on a new dataset.

Relying on a suitable pretraining, we are able to further enhance the training phase of a POD-DL-ROM, combining models of different fidelities (e.g., by considering coarser/finer spatial discretizations, as well as different physical laws, more/less parameters or larger/smaller parameter ranges).  In particular, we train the POD-DL-ROM neural network on an initial simpler task, obtaining a set of parameters $\boldsymbol{\theta}^*_S = \{ ({\bf W}_{S,i}^*, {\bf b}_{S,i}^*)\}_{i=1}^L$, and then use them to initialize the training of the POD-DL-ROM neural network on a more complex problem, by setting $\boldsymbol{\theta}^0_C= \boldsymbol{\theta}^*_S$. % (see \figurename~\ref{fig:Fig5-1}). 
Numerical results of Section \ref{sec:num_res} will show how pretraining, combined with the dimensionality reduction obtained through rPOD,  represents a cornerstone in view of drastically reducing the training computational time of a POD-DL-ROM, if compared to the time required for training, from scratch, a neural network on the more complex task.

%\begin{figure}[ht!]
%\centering
%\includegraphics[scale=0.45]{pretraining.pdf}
%\caption{Pretraining workflow.}
%\label{fig:Fig5-1}
%\end{figure}

\subsection{Extension to vector problems}

Compared to the DL-ROM technique, applied so far to scalar problems only, we have further generalized the POD-DL-ROM technique in order to handle vector problems, in analogy to what happens when treating  {\em red-green-blue} (RGB) images in general. DL algorithms. By considering the spatial discretization of a vector PDE problem, whose solution is a $d$-dimensional vector field,  problem (\ref{eq:FOM}) can be rewritten as
\begin{equation}
\label{eq:FOM_vectorial}
\begin{cases}
\mathbf{\dot{u}}_h^1(t;\boldsymbol{\mu}) = \mathbf{f}^1(t, \mathbf{u}_h^1(t;\boldsymbol{\mu}), \ldots, \mathbf{u}_h^d(t;\boldsymbol{\mu}); \boldsymbol{\mu}), \qquad t \in (0, T),\\
\vdots \\
\mathbf{\dot{u}}_h^d(t;\boldsymbol{\mu}) = \mathbf{f}^d(t, \mathbf{u}_h^1(t;\boldsymbol{\mu}), \ldots,\mathbf{u}_h^d(t;\boldsymbol{\mu}); \boldsymbol{\mu}), \qquad t \in (0, T),\\
\mathbf{u}_h^1(0;\boldsymbol{\mu})=\mathbf{u}_0^1(\boldsymbol{\mu}), \\
\vdots \\
\mathbf{u}_h^d(0;\boldsymbol{\mu})=\mathbf{u}_0^d(\boldsymbol{\mu}),
\end{cases}
\end{equation}
where $\mathbf{u}_h^i:[0,T) \times \mathcal{P} \rightarrow \mathbb{R}^{N_h^i}$ is the  solution of the $i$-th equation in (\ref{eq:FOM_vectorial}),  $\mathbf{u}_0^i : \mathcal{P} \rightarrow \mathbb{R}^{N_h^i}$ is the $i$-th initial datum, $\mathbf{f}^i : (0,T) \times  \mathbb{R}^{N_h^i}  \times \mathcal{P} \rightarrow \mathbb{R}^{N_h^i}$ describes the dynamics of $\mathbf{u}_h^i(t;\boldsymbol{\mu})$, $i= 1, \ldots, d$, with $d= 2,3$. Depending on the problem at hand, some of the equations appearing in (\ref{eq:FOM_vectorial}) might not involve the derivatives and related initial conditions; this is what happens, for instance, in the case of unsteady Navier-Stokes equations for incompressible flows, where the equation expressing flow incompressibility involves the velocity components, but no time derivatives.

Provided the solution of (\ref{eq:FOM_vectorial}), associated to a particular instance $(t, \boldsymbol{\mu})$, and the orthonormal basis $\mathbf{V}_{N,i} \in \mathbb{R}^{N_h^i \times N}$, with $i= 1, \ldots, d$, found though rSVD, we compute the intrinsic components $\mathbf{V}_{N,1}^T \mathbf{u}_h^1(t;\boldsymbol{\mu}), \ldots,$ $\mathbf{V}_{N,d}^T \mathbf{u}_h^d(t;\boldsymbol{\mu})$,
reshape each component in a square matrix of dimension $(\sqrt{N}, \sqrt{N})$, where $N = 2^{(2m)}$ with $m \in \mathbb{N}$, and stack them together forming a tensor with $d$ {\em channels}. \textcolor{black}{Thus, each vectorial component of the solution of problem (\ref{eq:FOM_vectorial}) is reshaped in a square matrix; then, they are stacked together  forming a tensor of dimension $(\sqrt{N}, \sqrt{N}, 3)$.}
% (a sketch in \figurename~\ref{fig:Fig1}).
%\begin{figure}[ht!]
%\centering
%\includegraphics[scale=0.55]{3D_image.pdf}
%\caption{Each vectorial component of the solution of problem (\ref{eq:FOM_vectorial}) is reshaped in a square matrix and they are stacked together thus forming a tensor of dimension $(\sqrt{N}, \sqrt{N}, 3)$.}
%\label{fig:Fig1}
%\end{figure}

This approach allows the dimensions $N_h^i$, $i= 1, \ldots, d$, of each FOM component, to be different. Indeed, it is the rPOD dimension $N$ used to reduce each vector component that must be kept equal for $i= 1, \ldots, d$. We remark that by stacking all components together allows to reduce the number of parameters, and then the training and testing computational times of POD-DL-ROM.

\subsection{Training and Testing Algorithms}

We summarize in this section both the training and the testing stage of the POD-DL-ROM technique. Regarding the setting of the optimization algorithm, the way to select the hyperparameters and the architecture of the neural networks, we refer to \cite{fresca2020comprehensive,fresca2020deep}. We denote by 
$\mathbf{M} \in \mathbb{R}^{(n_{\boldsymbol{\mu}} + 1) \times N_s}$ the matrix collecting all the parameter instances corresponding to the computed snapshots; these latter are included in the snapshot matrix ${\bf S}$, defined in \eqref{eq:snap_matrix1}. Note that, in the  case of a vector problem, the snapshot matrix $\mathbf{S} \in \mathbb{R}^{\sum_i N_h^i \times N_s}$ takes the form
\begin{equation*}
\mathbf{S} =
\begin{bmatrix}
\mathbf{S}_1 \\
\ldots \\
\mathbf{S}_d
\end{bmatrix}
=
\begin{bmatrix}
 \mathbf{u}_h^1(t^1; \boldsymbol \mu_1) \; | \; \ldots \; | \; \mathbf{u}_h^1(t^{N_t}; \boldsymbol \mu_1) \; | \; \ldots \; | \;
\mathbf{u}_h^1(t^1 ; \boldsymbol \mu_{N_{train}}) \; | \; \ldots \; | \; \mathbf{u}_h^1(t^{N_t} ; \boldsymbol \mu_{N_{train}}) \\
\ldots \\
 \mathbf{u}_h^d(t^1; \boldsymbol \mu_1) \; | \; \ldots \; | \; \mathbf{u}_h^d(t^{N_t}; \boldsymbol \mu_1) \; | \; \ldots \; | \;
\mathbf{u}_h^d(t^1 ; \boldsymbol \mu_{N_{train}}) \; | \; \ldots \; | \; \mathbf{u}_h^d(t^{N_t} ; \boldsymbol \mu_{N_{train}})
\end{bmatrix}.
\end{equation*}

The input and the output of the POD-DL-ROM are normalized by applying to each channel of the $d$-dimensional tensor $\mathbf{S}$ the affine transformation detailed below.  %Provided the training parameter matrix $\mathbf{M}^{train} \in \mathbb{R}^{(n_{\boldsymbol{\mu}} + 1) \times N_s}$
After splitting the data in $\mathbf{M} = [\mathbf{M}^{train}, \mathbf{M}^{val}]$ and $\mathbf{S} = [\mathbf{S}^{train}, \mathbf{S}^{val}]$ -- where $\mathbf{M}^{val}, \mathbf{S}^{val} \in \mathbb{R}^{\sum_i N_h^i \times \alpha N_s}$,  and $\alpha$ is a user-defined training-validation splitting fraction -- we define
\begin{equation}
\label{max_min_M}
M_{max}^i = \max_{j =1, \ldots, N_s} M_{ij}^{train},  \qquad %\textnormal{and} \quad 
M_{min}^i = \min_{j =1, \ldots, N_s} M_{ij}^{train},
\end{equation}
%where $M_{max}, M_{min} \in \mathbb{R}^{(n_{\boldsymbol{\mu}} + 1)}$, 
so that   parameters are normalized by applying the following transformation
\begin{equation}
\label{normalization}
M_{ij}^{train} \mapsto \frac{M_{ij}^{train} - M_{max}^i}{M_{max}^i - M_{min}^i}, \qquad
i = 1 \ldots, n_{\boldsymbol{\mu}} + 1,  \ j = 1, \ldots, N_s,
\end{equation} 
%for $$ and $$. 
that is, each feature of the training parameter matrix is rescaled according to its maximum and minimum values. Regarding instead the training snapshot matrix $\mathbf{S}^{train} \in \mathbb{R}^{\sum_i N_h^i \times N_s}$, we define
\begin{equation}
\label{max_min_S}
%\begin{split}
S_{max}^k = \max_{i =1, \ldots, N_h} \max_{j = 1, \ldots, N_s} S_{ij}^{train}, \quad \quad % \\ 
S_{min}^k = \min_{i =1, \ldots, N_h} \min_{j = 1, \ldots, N_s} S_{ij}^{train}, \quad \quad k = 1, \ldots, d
%\end{split}
\end{equation}
and apply transformation (\ref{normalization}),  by replacing $M_{max}^i, M_{min}^i$ with $S_{max}, S_{min} \in \mathbb{R}$ respectively, to each channel of $\mathbf{S}$ --  that is, we use the same maximum and minimum values for all the features of the snapshot matrix, as in \cite{carlberg2018model,gonzalez2018deep}. Using the latter approach or employing each feature's maximum and minimum values, for the matrix $\mathbf{S}^{train}$, does not lead to remarkable changes in the POD-DL-ROM performance. Transformation (\ref{normalization}) is applied also to the validation and testing sets, but considering as maximum and minimum the values computed over the training set. 
\textcolor{black}{In order to rescale the reconstructed solution to the original values, we apply the inverse transformation of (\ref{normalization}).} 

We detail the algorithms through which the training and the testing of the neural network are performed in Algorithms \ref{training_algorithm_POD} and \ref{testing_algorithm_POD}. 
During the training phase, the optimal parameters of the POD-DL-ROM are found by solving the optimization problem (\ref{eq:minimization_problem})-(\ref{eq:loss_N}) through the back-propagation and ADAM algorithms (see Algorithm \ref{training_algorithm_POD}). At testing time, the encoder function is instead discarded (see Algorithm \ref{testing_algorithm_POD}).  {By exploiting an early stopping criterion, we  stop the training if the loss function does not decrease over a certain number of epochs  over the validation set.}

\begin{algorithm}[t!]
\caption{POD-DL-ROM training algorithm}
\begin{algorithmic}[1]
\Require Parameter matrix $\mathbf{M} \in \mathbb{R}^{(n_{\boldsymbol{\mu}} + 1) \times N_s}$, snapshot matrix $\mathbf{S} \in \mathbb{R}^{\sum_i N_h^i \times N_s}$, training-validation splitting fraction $\alpha$, starting learning rate $\eta$, batch size $N_b$, maximum number of epochs $N_{epochs}$, number of minibatches $N_{mb} = (1 - \alpha)N_s/N_b$.
\Ensure  Optimal model parameters $\boldsymbol{\theta}^*= (\boldsymbol{\theta}_E^*, \boldsymbol{\theta}_{DF}^*, \boldsymbol{\theta}_D^*)$.
\vspace{0.3cm}
\State Compute rPOD basis matrix $\mathbf{V}_N = [\mathbf{V}_{N,1} | \ldots | \mathbf{V}_{N,d}]^T$ \;
\State Randomly shuffle $\mathbf{M}$ and $\mathbf{S}$ \;
\State Split data in $\mathbf{M} = [\mathbf{M}^{train}, \mathbf{M}^{val}]$ and $\mathbf{S} = [\mathbf{S}^{train}, \mathbf{S}^{val}]$ (with $\mathbf{M}^{val}, \mathbf{S}^{val} \in \mathbb{R}^{\sum_i N_h^i \times \alpha N_s}$)\;
\State Compute intrinsic coordinates $\mathbf{S}^{train}_N = [\mathbf{S}^{train}_1| \ldots | \mathbf{S}^{train}_d]^T$ where $\mathbf{S}^{train}_{i} = \mathbf{V}_{N,i}^T \mathbf{S}^{train}_i$, $i = 1, \ldots, d$ \;
\State Compute intrinsic coordinates $\mathbf{S}^{val}_N = [\mathbf{S}^{val}_1| \ldots | \mathbf{S}^{val}_d]^T$ where $\mathbf{S}^{val}_{i} = \mathbf{V}_{N,i}^T \mathbf{S}^{val}_i$, $i = 1, \ldots, d$ \;
\State Normalize data in $\mathbf{M}$ and $\mathbf{S}_N = [\mathbf{S}^{train}_N, \mathbf{S}^{val}_N]$\;
\State Randomly initialize $\boldsymbol{\theta}^0=(\boldsymbol{\theta}_{E}^0, \boldsymbol{\theta}_{DF}^0, \boldsymbol{\theta}_{D}^0)$\;
\State $n_{e} = 0$
\While{($\neg$early-stopping \textbf{and} $n_{e} \le N_{epochs}$)}
\For{$k = 1 : N_{mb}$}
    \State Sample a minibatch $(\mathbf{M}^{batch}, \mathbf{S}_N^{batch}) \subseteq (\mathbf{M}^{train}, \mathbf{S}^{train}_N)$\;
    \State $\mathbf{S}^{batch}_N =$ reshape$(\mathbf{S}^{batch}_N) \in \mathbb{R}^{N_b \times \sqrt{N} \times \sqrt{N} \times d}$ \;
    \State $\mathbf{\widetilde{S}}^{batch}_n(\boldsymbol{\theta}_{E}^{N_{mb} n_{e} + k}) = \mathbf{f}_{n}^E(\mathbf{S}^{batch}_N; \boldsymbol{\theta}_{E}^{N_{mb} n_{e} + k})$\;
    \State $\mathbf{S}^{batch}_n(\boldsymbol{\theta}_{DF}^{N_{mb} n_{e} + k}) = \boldsymbol{\phi}_n^{DF}(\mathbf{M}^{batch}; \boldsymbol{\theta}_{DF}^{N_{mb} n_{e} + k})$\;
    \State $\mathbf{\widetilde{S}}^{batch}_N(\boldsymbol{\theta}_{DF}^{N_{mb} n_{e} + k}, \boldsymbol{\theta}_{D}^{N_{mb} n_{e} + k}) = \mathbf{f}_{N}^D(\mathbf{S}^{batch}_n(\boldsymbol{\theta}_{DF}^{N_{mb} n_{e} + k}); \boldsymbol{\theta}_{D}^{N_{mb} n_{e} + k})$
    \State $\mathbf{\widetilde{S}}^{batch}_N =$ reshape$(\mathbf{\widetilde{S}}^{batch}_N) \in \mathbb{R}^{N_b  \times N \times d}$\;
    \State Accumulate loss (\ref{eq:loss_N})  {on $(\mathbf{M}^{batch}, \mathbf{S}_N^{batch})$} and compute $\widehat{\nabla}_{\theta} \mathcal{J}$\;
  	\State $\boldsymbol{\theta}^{N_{mb} n_{e} + k + 1} = \textnormal{ADAM}(\eta, \widehat{\nabla}_{\theta} \mathcal{J}, \boldsymbol{\theta}^{N_{mb} n_{e} + k})$\;
  \EndFor
%  \State Repeat lines 8-12 $(M^{val}, S^{val})$  {considering updated weights}
  \State Repeat instructions 12-16 on $(\mathbf{M}^{val}, \mathbf{S}^{val}_N)$ with the updated weights $\boldsymbol{\theta}^{N_{mb} n_{e} + k + 1}$
  \State Accumulate loss (\ref{eq:loss_N}) on $(\mathbf{M}^{val}, \mathbf{S}^{val}_N)$ to evaluate early-stopping criterion
  \State $n_{e} = n_{e} + 1$
\EndWhile
\end{algorithmic}
\label{training_algorithm_POD}
\end{algorithm}
\begin{algorithm}[ht!]
\caption{POD-DL-ROM testing algorithm}
\begin{algorithmic}[1]
\Require Testing parameter matrix $\mathbf{M}^{test} \in \mathbb{R}^{(n_{\boldsymbol{\mu}} + 1) \times (N_{test} N_t)}$, rPOD basis matrix $\mathbf{V}_N$, % optimal model parameters 
$(\boldsymbol{\theta}_{DF}^*, \boldsymbol{\theta}_D^*)$.
\Ensure ROM approximation matrix $\mathbf{\widetilde{S}}_h \in \mathbb{R}^{\sum_i N_h^i \times (N_{test} N_t)}$.
\vspace{0.3cm}
\State Load $\boldsymbol{\theta}_{DF}^*$ and $\boldsymbol{\theta}_D^*$\;
\State $\mathbf{S}_n(\boldsymbol{\theta}_{DF}^*) = \boldsymbol{\phi}_n^{DF}(\mathbf{M}^{test}; \boldsymbol{\theta}_{DF}^*)$\;
\State $\mathbf{\widetilde{S}}_N(\boldsymbol{\theta}_{DF}^*, \boldsymbol{\theta}_{D}^*) = \mathbf{f}_{N}^D(\mathbf{S}_n(\boldsymbol{\theta}_{DF}^*); \boldsymbol{\theta}_{D}^*)$
\State $\mathbf{\widetilde{S}}_N=$ reshape$(\mathbf{\widetilde{S}}_N)$\;
\State $\mathbf{\widetilde{S}}_h= \mathbf{V}_N\mathbf{\widetilde{S}}_N$
\end{algorithmic}
\label{testing_algorithm_POD}
\end{algorithm}

\textcolor{black}{We remark that with $\mathbf{\widetilde{S}}_n$ we refer to a matrix collecting in its columns the output of the encoder function of the convolutional AE applied to each column of the snapshot matrix $\mathbf{S}$. In the same way,  the columns of $\mathbf{{S}}_n$  collect the minimal coordinates, output of the DFNN, for each sample in the parameter matrix $\mathbf{M}$, and $\mathbf{\widetilde{S}}_h$ is a matrix whose columns are the intrinsic coordinates approximations, outputs of the decoder function of the convolutional AE, associated to the columns of  $\mathbf{{S}}_n$}.

\section{Numerical results}
\label{sec:num_res}

We assess the numerical performance of the proposed POD-DL-ROM technique, by focusing on the training and testing computational times required to construct and deploy a POD-DL-ROM, and the use of pretraining, on four different linear or nonlinear parametrized PDE problems: {\em (i)} a linear unsteady advection-diffusion-reaction equation; {\em (ii)} the monodomain system for cardiac electrophysiology; {\em (iii)} a nonlinear elastodynamics problem, and {\em (iv)} the unsteady Navier-Stokes equations for incompressible flows. \\ 

 To evaluate the performance of the POD-DL-ROM technique, we rely on two error indicators:
\begin{itemize}
\item the error indicator $\epsilon_{rel} \in \mathbb{R}$ given by
\begin{equation}
\epsilon_{rel}(\mathbf{u}_h, \mathbf{\tilde{u}}_h) = \frac{1}{N_{test}} \sum_{i  = 1}^{N_{test}} \left(\displaystyle \frac{\sqrt{ \sum_{k=1}^{N_t} || \mathbf{u}^k_h(\boldsymbol{\mu}_{test,i}) - \mathbf{\tilde{u}}^k_h(\boldsymbol{\mu}_{test,i}) ||^2}}{\sqrt{\sum_{k=1}^{N_t} || \mathbf{u}_h^k(\boldsymbol{\mu}_{test,i}) ||^2}} \right),
\label{eq:error_indicator}
\end{equation}
\item the relative error $\boldsymbol{\epsilon}_k \in \mathbb{R}^{\sum_i N_h^i}$, for $k = 1, \ldots, N_t$, defined as
\begin{equation}
\displaystyle \boldsymbol{\epsilon}_k(\mathbf{u}_h, \mathbf{\tilde{u}}_h) = \displaystyle \frac{ | \mathbf{u}^k_h(\boldsymbol{\mu}_{test}) - \mathbf{\tilde{u}}^k_h(\boldsymbol{\mu}_{test}) |}{\sqrt{\frac{1}{N_t}\sum_{k=1}^{N_t} || \mathbf{u}^k_h(\boldsymbol{\mu}_{test}) ||^2}}.
\label{eq:relative_error}
\end{equation}
\end{itemize}
The coefficient $\omega_h \in [0,1]$ in (\ref{eq:loss_N}) is set equal to 0.5 according to  the results shown in \cite{fresca2020comprehensive}, where a detailed analysis suggested to select values of $\omega_h$ equidistant from the extrema of $[0,1]$.  The rPOD dimension $N$ is selected, in all test cases, in such a way that $\epsilon_{rel}(\mathbf{u}_h, \mathbf{V}_N\mathbf{V}_N^T \mathbf{u}_h) \approx 10^{-4}$, whereas the dimension of the nonlinear trial manifold $n$ is set trying to match the dimension of the solution manifold $n_{\mu} + 1$. The POD-DL-ROM neural network is implemented by means of the \texttt{Tensorflow} DL framework \cite{abadi2016tensorflow}.

\subsection{Test 1: unsteady advection-diffusion-reaction equation}

The first test case we consider deals with the solution $u=u({\bf x}, t; \boldsymbol{\mu})$ of the following advection-diffusion-reaction system
\begin{equation}
\left\{
\begin{aligned}
& \frac{\partial u}{\partial t}  - \textnormal{div}( \mu_1  \nabla u) + \mathbf{b}(t;\mu_2) \cdot \nabla u + c u   =  f(\mu_3, \mu_4) & \  & (\mathbf{x}, t) \in \Omega \times (0,T),\\
& \mu_1 \nabla u  \cdot \mathbf{n} = 0 &\ & (\mathbf{x}, t) \in \partial \Omega \times (0,T), \\
& u(0) = 0 & \  &  \mathbf{x} \in \Omega,
\end{aligned}
\right.
\label{eq:ADR}
\end{equation}
in the two-dimensional domain $\Omega = (0, 1)^2$, where
\begin{equation*}
f(\mathbf{x}; \mu_3, \mu_4) = 10 \exp(-((x-\mu_3)^2 + (y-\mu_4)^2)/0.07^2)
\end{equation*}
and
\begin{equation*}
\mathbf{b}(t; \mu_2) = [\cos(\pi/\mu_2t), \sin(\pi/\mu_2t)]^T.
\end{equation*}
We consider $n_{\mu}=4$ parameters, belonging to  $\mathcal{P} = [0.002, 0.005] \times [30, 70] \times [0.4, 0.6]^2$; we build a FOM considering a space discretization made by linear $(\mathbb{P}_1)$ finite elements, considering $N_h = 10657$ DOFs, and  a Backward Differentiation Formula (BDF) of order 2 considering a time step  $\Delta t = 2\pi/20$ over  $(0,T)$ with $T=10\pi$, as time discretization. For different values of  $\mu_3$ and $\mu_4$, the solution of (\ref{eq:ADR}) exhibits different patterns, due to the location of the distributed source; the dependence on $\mu_1$ and $\mu_2$ impact instead on the relative importance of diffusion and advection terms, and on the direction of this latter. We expect that, for the case at hand, POD-Galerkin ROMs might involve a large number of basis functions, also because; note also that the dependence of the solution on $\mu_3$ and $\mu_4$ is nonlinear (that is, the problem is nonaffinely parametrized).

Regarding the construction of the proposed POD-DL-ROM, for the training of the neural networks, we consider   $N_t = 100$ time instances   and $N_{train} = 5 \times 5 \times 5 \times 4 = 500$ training-parameter instances, uniformly distributed in each parametric direction. At testing phase,  $N_{test} = 4 \times 4 \times 4 \times 3 = 192$ testing-parameter instances have been considered instead, different from the training ones. The maximum number of epochs is set equal to $N_{epochs} = 10000$, the batch size is $N_b = 120$ and, regarding the early-stopping criterion, we stop the training if the loss function does not decrease within 500 epochs.    {We set $N = 64$ as dimension of the rPOD basis (i.e. the linear trial manifold generated by means of rSVD), and $n = n_{\mu} + 1 =5$ as  dimension of the reduced nonlinear trial manifold.} The training and testing phases of the POD-DL-ROM neural network have been performed on a Tesla V100 32GB GPU.

\begin{figure}[ht!]
\centering
\includegraphics[scale=0.27]{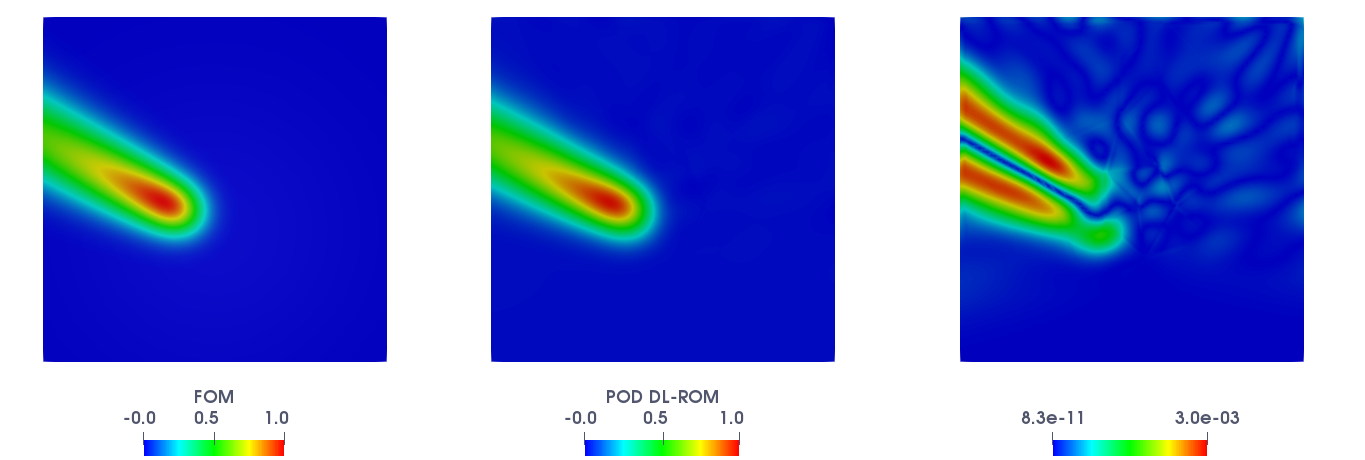} \\
\vspace{0.05cm}
\includegraphics[scale=0.27]{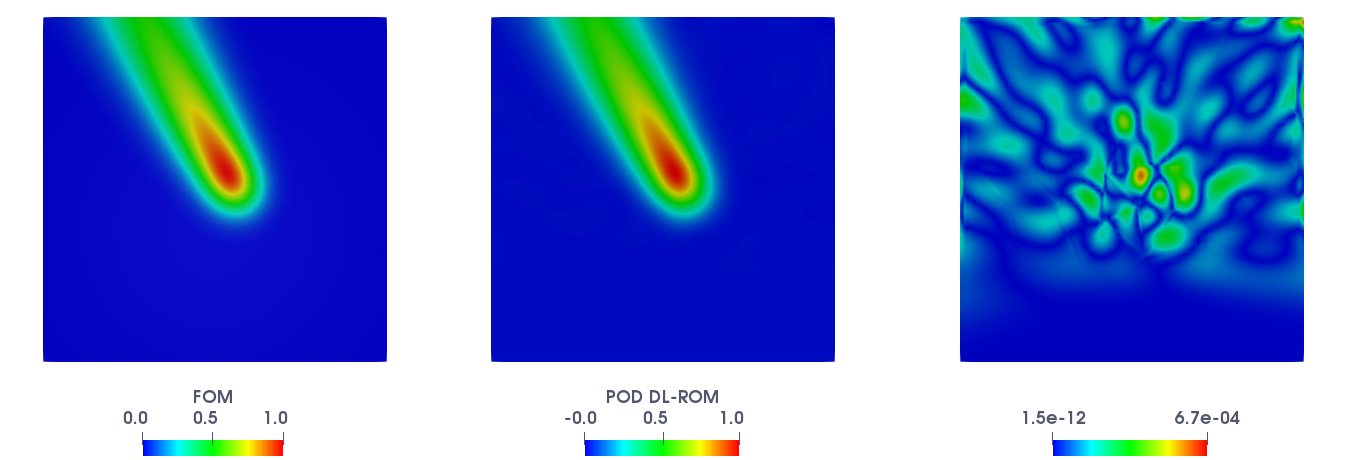}
\vspace{-0.15cm}
\caption[\textit{Test 1}: FOM, POD-DL-ROM, with $n=5$ and $N=64$, solutions and relative error $\boldsymbol{\epsilon}_k$, for the testing-parameter instances $\boldsymbol{\mu}_{test} = (0.425, 0.425, 35, 0.0045)$ and $\boldsymbol{\mu}_{test} = (0.575, 0.475, 45, 0.0045)$ at $t = 29.53$.]{\textit{Test 1}: FOM (left), POD-DL-ROM (center), with $n=5$ and $N=64$, solutions and relative error $\boldsymbol{\epsilon}_k$ (right), for the testing-parameter instances $\boldsymbol{\mu}_{test} = (0.425, 0.425, 35, 0.0045)$ (top) and $\boldsymbol{\mu}_{test} = (0.575, 0.475, 45, 0.0045)$ (bottom) at $t = 29.53$.}
\label{fig:Fig5-2}
\end{figure}

In \figurename~\ref{fig:Fig5-2} we show the FOM  and the POD-DL-ROM solutions, for the testing-parameter instances $\boldsymbol{\mu}_{test}=(0.425, 0.425, 35, 0.0045)$  and $\boldsymbol{\mu}_{test}=(0.575, 0.475, 45, 0.0045)$ at $t = 29.53$, respectively, together with the relative error (\ref{eq:relative_error}).

\begin{figure}[b!]
\vspace{-0.25cm}
\centering
\includegraphics[scale=0.25]{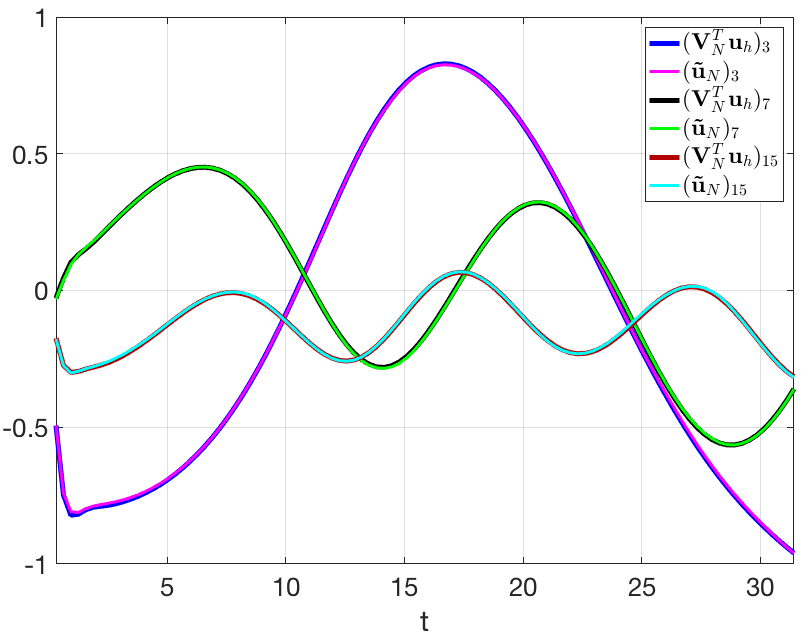} \hspace{0.5cm}
\includegraphics[scale=0.25]{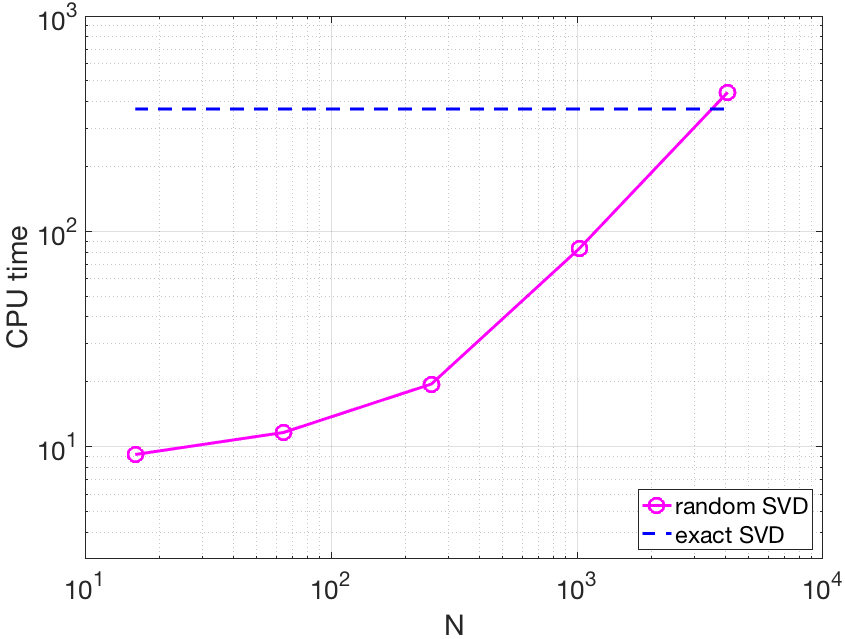}
\vspace{-0.25cm}
\caption{\textit{Test 1}. Left: comparison between the intrinsic coordinates $\mathbf{V}_N^T\mathbf{u}_h$ components and the POD-DL-ROM approximation $\mathbf{\tilde{u}}_N$ for the testing-parameter instance $\boldsymbol{\mu}_{test} = (0.575, 0.475, 45, 0.0045)$. Right:  SVD and rSVD CPU times vs. $N$.}
\label{fig:Fig5-3}
\end{figure}

The comparison between some components of the intrinsic coordinates vector $\mathbf{V}_N^T\mathbf{u}_h (t; \boldsymbol{\mu}_{test})$ and their POD-DL-ROM approximation, for the testing parameter instance $\boldsymbol{\mu}_{test} = (0.575, 0.475,$ $45, 0.0045)$, is shown is \figurename~\ref{fig:Fig5-3}. We remark that, as expected, the first components are the ones retaining most of the energy of the system; thus being the ones with higher magnitude \cite{quarteroni2016reduced}.

In \figurename~\ref{fig:Fig5-3} (right) we show the CPU  times, as function of $N$, required by (exact) SVD and  rSVD to compute the linear POD space, and thus performing the first level of dimensionality reduction required by the POD-DL-ROM. The CPU time required by SVD is not affected by $N$, while the one required  by rSVD increases with respect to $N$ -- note that the two times  almost coincide for $N = 4096$, a dimension for which constructing a ROM could in principle be avoided. Hence, rSVD is always preferable with respect to SVD; in particular, for the choice $N = 64$, using rSVD allows a speed-up equal to 32 with respect to SVD.

\begin{figure}[t!]
\centering
\hspace{0.3cm}
\includegraphics[scale=0.14]{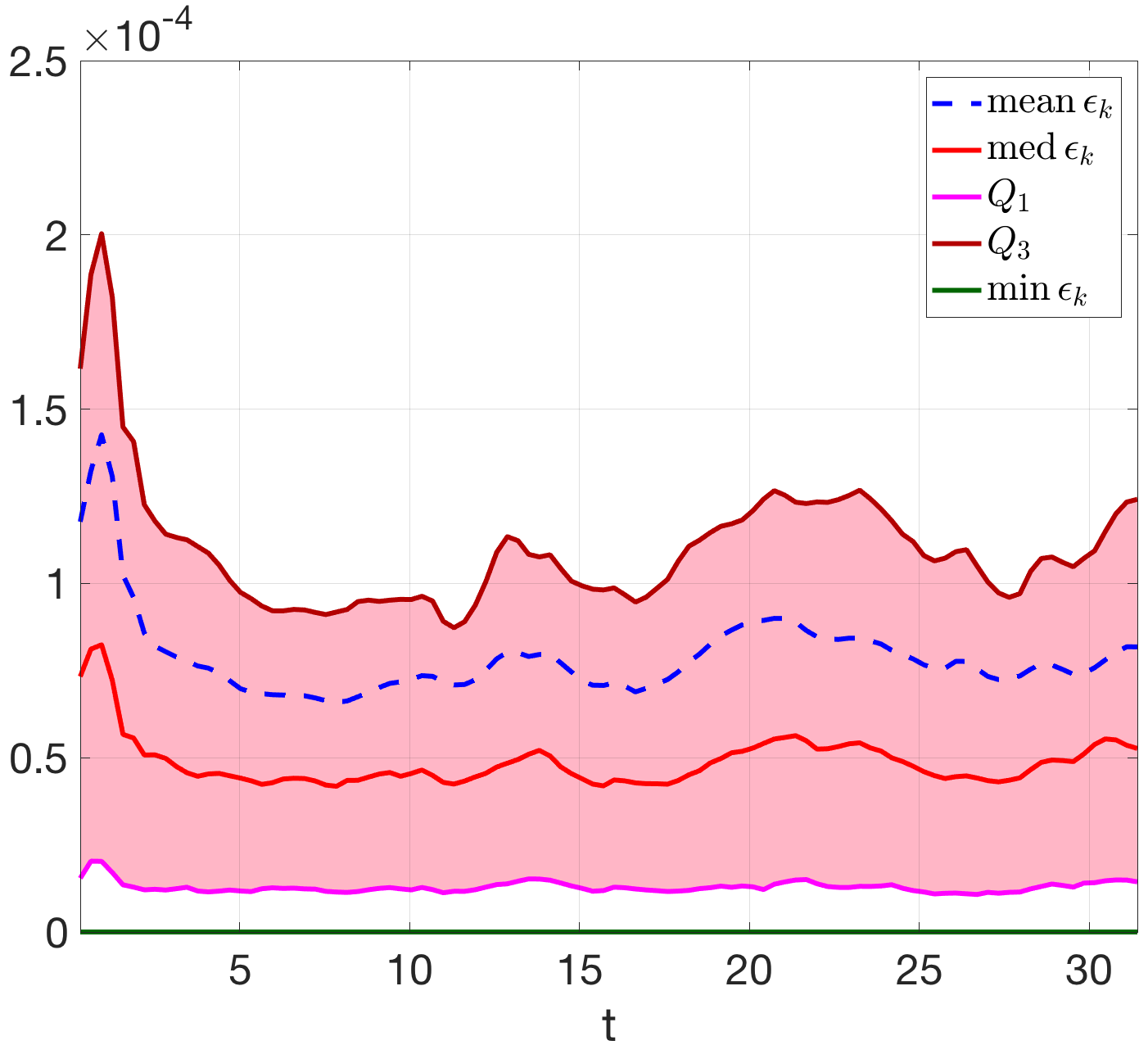} \hspace{0.85cm}
\includegraphics[scale=0.255]{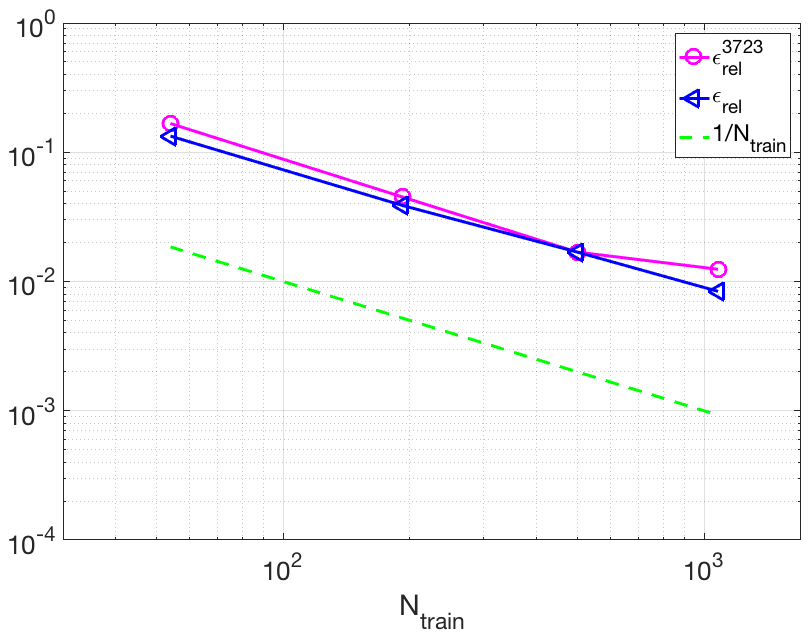}
%\caption[\textit{Test 1}. Left: SVD and rSVD CPU times vs. $N$ and trend of the relative error over time for $\boldsymbol{\mu}_{test} = ( 0.575, 0.475, 45, 0.0045 )$.]
\caption{\textit{Test 1}. Left: trend of the relative error over time for $\boldsymbol{\mu}_{test} = (0.575, 0.475, 45, 0.0045)$. Right: error indicator $\epsilon_{rel}$ vs. $N_{train}$.} %]{\textit{Test 1}: Error indicator $\epsilon_{rel}$ vs. $N_{train}$.}}
\label{fig:Fig5-4}
\end{figure}

The trend of the relative error (\ref{eq:relative_error}) over time, for the selected testing-parameter instance $\boldsymbol{\mu}_{test} = (0.575, 0.475, 45, 0.0045)$, is displayed in \figurename~\ref{fig:Fig5-4} (left), where the mean (over the domain), the median,  the first and third quartile of the relative error, as well as its minimum, are reported. The interquartile range (IQR) shows that the distribution of the error is almost uniform over time.

\textcolor{black}{We also analyze the convergence properties of the POD-DL-ROM  by varying the number of training-parameter instances provided to the neural network. In particular, in \figurename~\ref{fig:Fig5-4} (right) we report the trend of the error indicator (\ref{eq:error_indicator}), over the testing set, versus $N_{train}$, i.e. the size of the training dataset. With $\epsilon_{rel}^{3723}$ we refer to the value of the error indicator obtained by setting the maximum number of epochs equal to $N_{epochs}=3723$, which are the iterations performed during the training phase by considering $N_{train}=500$. By increasing $N_{train}$ , that is by providing more data to the POD DL-ROM neural network, its approximation capability increases, thus yielding  a decrease in the error indicator. In particular, the  error indicator  (\ref{eq:relative_error}) decays with a rate that is  about $1/N_{train}$.}

\begin{figure}[ht!]
\vspace{-0.15cm}
\centering
\includegraphics[scale=0.25]{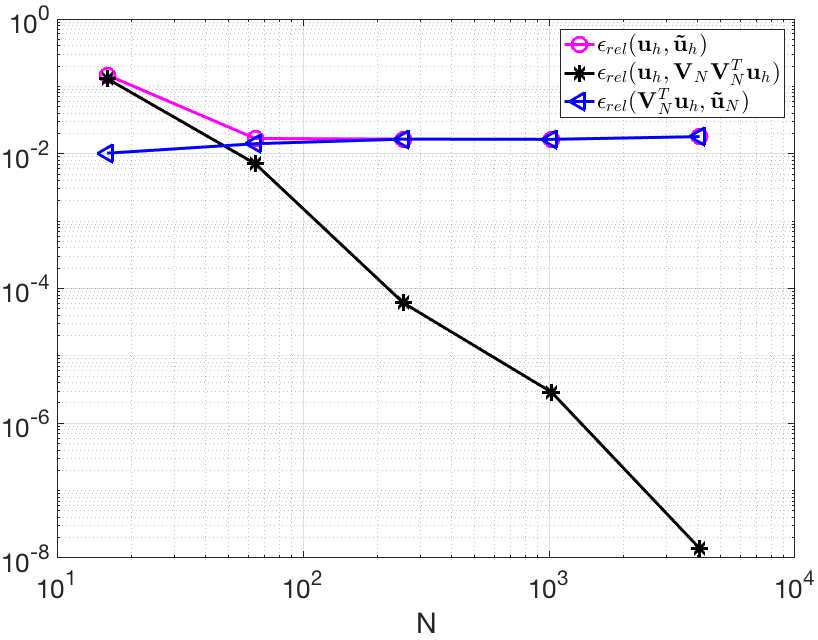}
\hspace{0.3cm}
\includegraphics[scale=0.25]{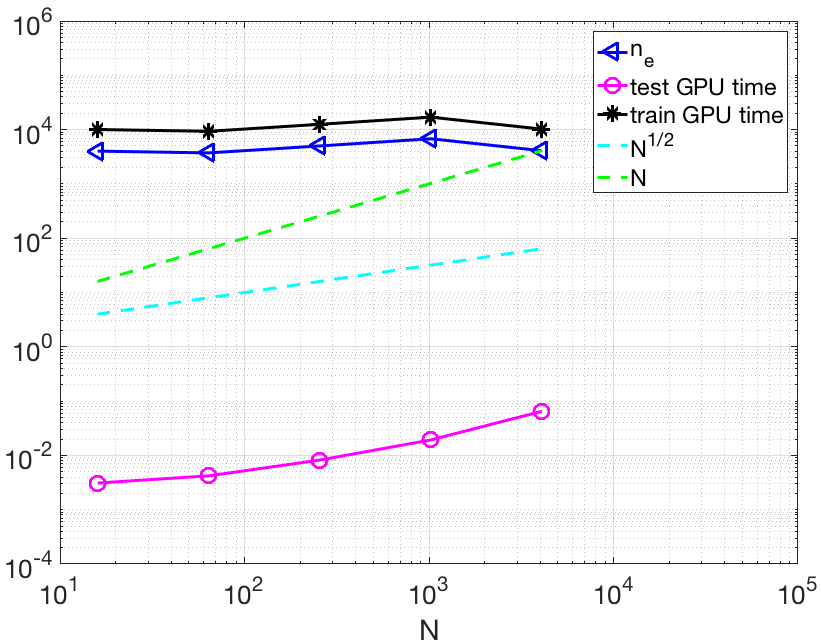}
\vspace{-0.25cm}
\caption[\textit{Test 1}: Error indicator $\epsilon_{rel}$ vs. $N$ and GPU training and testing computational times vs. $N$.]{\textit{Test 1}: Left: Error indicator $\epsilon_{rel}$ vs. $N$. Right: GPU training and testing computational times vs. $N$.}
\label{fig:Fig5-5}
\end{figure}

We finally assess the accuracy and the efficiency of the POD-DL-ROM with respect to the dimension $N$ of the POD space. In \figurename~\ref{fig:Fig5-5} (left) we display the error indicator (\ref{eq:error_indicator}) computed on the FOM and POD-DL-ROM solutions, on the FOM and the optimal-POD (i.e., the projection of the FOM solution onto the linear trial manifold generated through rSVD) solutions,   the intrinsic coordinates $\mathbf{V}_N^T \mathbf{u}_h$, and the approximated ones $\mathbf{\tilde{u}}_N$. The trend of $\epsilon_{rel}(\mathbf{u}_h, \mathbf{\tilde{u}}_h)$, for $N = 16$, is dictated by the projection error $\epsilon_{rel}(\mathbf{u}_h, \mathbf{V}_N\mathbf{V}_N^T\mathbf{u}_h)$, thus indicating that, for $N=16$, the rPOD dimension is too small to accurately reconstruct the FOM solution. For $N \ge 64$, the error indicator $\epsilon_{rel}(\mathbf{u}_h, \mathbf{\tilde{u}}_h)$ remains almost  unchanged, and almost coincide with $\epsilon_{rel}(\mathbf{V}_N^T\mathbf{u}_h, \mathbf{\tilde{u}}_N)$. As observed for the DL-ROM \cite{fresca2020comprehensive}, this behavior is related to the fact that an increase of $N$ only implies the addition of few parameters to the POD-DL-ROM neural network, i.e. the approximation capability of the network remains almost the same while the input increases, thus resulting in a more difficult task.

%\begin{figure}[ht!]
%\centering
%\includegraphics[scale=0.25]{ADR_err_vs_Ntrain.png}
%\caption[\textit{Test 1}: Error indicator $\epsilon_{rel}$ vs. $N_{train}$]{\textit{Test 1}: Error indicator $\epsilon_{rel}$ vs. $N_{train}$.}
%\label{fig:Fig5-5_0}
%\end{figure}

The GPU training and testing computational times versus $N$ are pointed out in \figurename~\ref{fig:Fig5-5} (right). The training time refers to the total time required for the training and validation phases; for the sake of completeness we also show the number of epochs $n_e$ along with $N$. The training time varies between 2 h 30 m and 4 h 40 m, we remark that we deal with $N_{train}=500$ training-parameter instances. The testing time consists instead in the time required to compute $N_t$ time instances for a testing-parameter instance. The trend is proportional to $N^{1/2}$ and, for example, for $N=64$ the testing time is equal to 4.2 $\times 10^{-3}$ s \textcolor{black}{thus leading to a speed-up $1.2 \times 10^{4}$ with respect to the solution of the FOM on a MacBook Pro Intel Core i7 6-core with 16 GB RAM.}

\subsection{Test 2: coupled PDE-ODE Monodomain/Aliev-Panfilov system}

We now consider a coupled PDE-ODE nonlinear system modeling the electrical behavior of the cardiac tissue, from the cellular scale to the tissue level: the Monodomain equation \cite{franzone2014mathematical} coupled with the Aliev-Panfilov ionic model \cite{AlievPanfilov}, in a square slab of  tissue $\Omega = (0, 10)^2$ cm:
\begin{equation}
\label{eq:monodomain}
\begin{cases}
\displaystyle \frac{\partial u}{\partial t} - \textnormal{div}({\bf D} \nabla u) +Ku(u - a)(u - 1) + uw = I_{app}(\mathbf{x}, t) \quad & (\mathbf{x}, t) \in \Omega \times (0,T), \vspace{0.1cm} \\
\displaystyle \frac{\partial w}{\partial t} +\displaystyle \Big( \epsilon_0 + \frac{c_1w}{c_2+u} \Big)(- w - Ku(u - b -1))=0 \quad & (\mathbf{x}, t) \in \Omega \times (0,T), \vspace{0.1cm}\\
\displaystyle \nabla u \cdot \mathbf{n} = 0 \quad & (\mathbf{x}, t) \in \partial \Omega \times (0,T),\\
u(\mathbf{x},0)=0, \; w(\mathbf{x},0)= 0 \quad & \mathbf{x} \in \Omega.
\end{cases}
\end{equation}
We consider  two $(n_{\mu}=2)$ parameters, consisting in the electric conductivities in the longitudinal and the transversal directions to the fibers, i.e., the conductivity tensor $\mathbf{D}(\mathbf{x}; \boldsymbol{\mu}) $ takes the form
\begin{equation}
\label{conductivity_tensor}
\mathbf{D}(\mathbf{x}; \boldsymbol{\mu}) = \mu_2 I + (\mu_1 - \mu_2) \mathbf{f}_0(\mathbf{x}) \otimes \mathbf{f}_0(\mathbf{x}),
\end{equation}
where $\mathbf{f}_0=(1, 0)^T$ and the parameter space is   $\mathcal{P} = 12.9 \cdot [0.06, 0.2]  \times 12.9 \cdot [0.03,0.1 ]$ $\textnormal{cm\textsuperscript{2}/ms}$. The applied current is defined as
\begin{equation*}
I_{app}(\mathbf{x}, t) = \frac{C}{2 \pi \alpha} \exp \bigg( - \frac{||\mathbf{x}||^2}{2 \beta} \bigg)\mathbf{1}_{[0, t]}(t),
\end{equation*}
where $C = 100$ mA, $\alpha = 1$, $\beta = 1$ cm\textsuperscript{2} and $t = 2$ ms. The parameters of the Aliev-Panfilov ionic model are set to $K = 8$, $a = 0.01$, $b = 0.15$, $\varepsilon_0 = 0.002$, $c_1 = 0.2$, and $c_2 = 0.3$, see, e.g., \cite{goktepe2010atrial}. The equations have been discretized in space through linear  $(\mathbb{P}_1)$ finite elements   by considering $N_h = 64 \times 64 = 4096$ grid points. For the time discretization and the treatment of nonlinear terms, we use a one-step, semi-implicit, first order scheme (see, e.g., \cite{pagani2018numerical} for further details) by considering a time step  $\Delta t = 0.1$ ms over the interval $(0,T)$, with $T=400$ ms.

For the training phase, we uniformly sample $N_t = 1000$ time instances in the interval $(0,T)$ and consider $N_{train} = 25$ training-parameters, i.e. $\boldsymbol{\mu}_{train}=12.9 \cdot (0.06 + i 0.035, 0.03 + j 0.0175)$ with $i,j = 0, \ldots, 4$. For the testing phase,  $N_{test} = 16$ testing-parameter instances have been considered, each of them given by $\boldsymbol{\mu}_{test}=12.9 \cdot (0.0775 + i 0.035, 0.0387 + j 0.0175)$ with $i,j = 0, \ldots, 3$. The maximum number of epochs is  $N_{epochs} = 10000$, the batch size is $N_b = 40$ and, regarding the early-stopping criterion, we stop the training if the loss function does not decrease \textcolor{black}{within} 500 epochs. The simulations are performed on a GTX 1070 8GB GPU. We considered problem (\ref{eq:monodomain}) first in \cite{fresca2020comprehensive} and we now compare the increased efficiency entailed by the use of POD-DL-ROM compared to a DL-ROM and to a POD-Galerkin ROM.

\begin{figure}[ht!]
\centering
\includegraphics[scale=0.225]{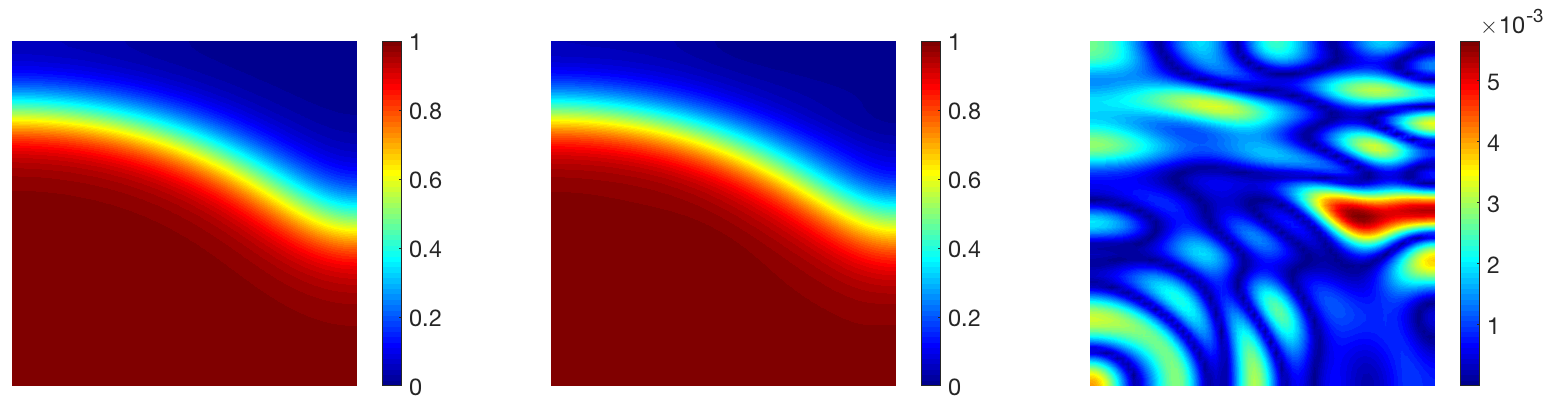}
\includegraphics[scale=0.225]{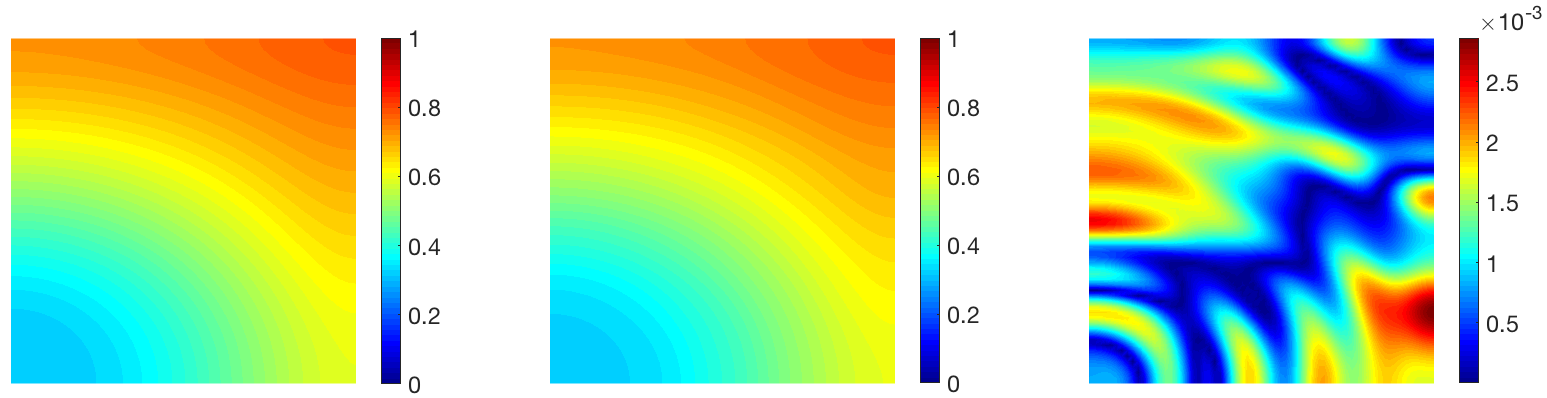}
\vspace{-0.25cm}
\caption[\textit{Test 2}: FOM, POD-DL-ROM solutions and relative error $\boldsymbol{\epsilon}_k$  for the testing-parameter instance $\boldsymbol{\mu}_{test} = 12.9 \cdot (0.1825, 0.0912)$ $\textnormal{cm\textsuperscript{2}/ms}$ at $t = 47.7$ ms and $t = 379.7$ ms, with $n=3$.]{\textit{Test 2}: FOM (left), POD-DL-ROM (center) solutions and relative error $\boldsymbol{\epsilon}_k$ (right) for the testing-parameter instance $\boldsymbol{\mu}_{test} = 12.9 \cdot (0.1825, 0.0912)$ $\textnormal{cm\textsuperscript{2}/ms}$ at $t = 47.7$ ms and $t = 379.7$ ms, with $n=3$.}
\label{fig:Fig6-6}
\end{figure}

In \figurename~\ref{fig:Fig6-6} we report the FOM and POD-DL-ROM solutions, the latter obtained by setting $n = 3$ and $N = 64$, along with the relative error (\ref{eq:relative_error}), for the testing-parameter instance $\boldsymbol{\mu}_{test} = 12.9 \cdot (0.1825, 0.0912)$ $\textnormal{cm\textsuperscript{2}/ms}$ at $t = 47.7$ ms (top) and $t = 379.7$ ms (bottom).

\textcolor{black}{We point out that setting $\omega_h$ in (\ref{eq:loss_N}) equal to 0.5, that is, performing a second level of dimensionality reduction by matching the dimension of the minimal components of the problem, leads to higher accuracies with respect to the case $\omega_h = 1$. In this respect, in \tablename~\ref{tb:tb0}, we show the values of the error indicator $\epsilon_{rel}$, over the testing set, by solving the optimization problem (\ref{eq:minimization_problem}) with and without considering the  term $\mathcal{L}_{int}(t^k, \boldsymbol{\mu}_i;  {\boldsymbol{\theta}})$ in (\ref{eq:loss_N}). In particular, setting $\omega_h=0.5$ allows to halve the error.
\begin{table}[h!]
\centering
\begin{tabular}{|c|c|c|}
\hline
 & $\omega_h=0.5$ & $\omega_h=1$ \\ \hline
$\epsilon_{rel}$ & $4.03 	\times 10^{-3}$ & $7.69 	\times 10^{-3}$ \\ \hline
\end{tabular}
\caption{\textit{Test 2}: Error indicator $\epsilon_{rel}$ for $\omega_h = 0.5$ and 1.}
\label{tb:tb0}
\end{table}}

The comparison among the DL-ROM and POD-DL-ROM neural networks number of parameters, the GPU training and validation time for one epoch, the total number of epochs, the total GPU training and validation (total time) and testing GPU computational times, \textcolor{black}{and the speed-up obtained, at testing time, with respect to the solution of the FOM\footnote{The FOM simulation is carried out on a MacBook Pro Intel Core i7 6-core with 16 GB RAM.}} are reported in \tablename~\ref{tb:tb1}. We also report the total CPU offline and online computational times\footnote{Here we employ a full 64 GB node (20 Intel\textsuperscript{\textregistered} Xeon\textsuperscript{\textregistered} E5-2640 v4 2.4GHz cores) of the HPC cluster available at MOX, Politecnico di Milano.}  required by the POD-Galerkin ROM with $N_c =6$ clusters, which corresponds to the choice providing the most efficient results (see, e.g., Test 4 of \cite{fresca2020comprehensive}), by keeping for all the models the same degree of accuracy $\epsilon_{rel} = 4.03 \times 10^{-3}$ and  running the code on the hardware it is optimized for. 

The use of a POD-DL-ROM not only entails even faster testing computational times with respect to a DL-ROM, due to the remarkable reduction of the number of parameters of the neural network, \textcolor{black}{but also reduces the training time of a factor 37.5 (resp. 5) compared to a  POD-Galerkin ROM (resp. a DL-ROM)}. The POD-DL-ROM thus results to be the most efficient ROM both at training and testing stages.
\begin{table}[h!]
\centering
\begin{tabular}{|l|l|l|l|l|l|l|}
\hline
 & $\#$params & \makecell{ train - val \\ $[s/epoch]$ } & $\#$epochs & total time & test $[s]$ & speed-up\\ \hline
DL-ROM (GPU) & 2342595 & 7.5 - 0.8 & 6981 & 15 h & 0.08 & $3.03 \times 10^3$\\ \hline
POD-DL-ROM (GPU) & 269057 & 1.5 - 0.15 & 866 & 24 m & 0.015 & $1.62 \times 10^{4}$\\ \hline
POD-Galerkin ROM ($N_c=6$) & - & -  & - & 115 m & 8 & $3.04 \times 10^1$\\ \hline
\end{tabular}
\caption{\textit{Test 2}: DL-ROM, POD-DL-ROM and POD-Galerkin ROM computational times.}
\label{tb:tb1}
\end{table}

We now investigate the use of pretraining, introduced in Section \ref{sec:pretraining}, in two different scenarios:
\begin{itemize}
\item when increasing the FOM dimension $N_h$;
\item when increasing the dimension of the parameter space $\mathcal{P}$.
\end{itemize}
First, we use the optimal parameters, weights and biases, of the POD-DL-ROM neural network found in the case of a FOM dimension $N_h = 4096$, to initialize the POD-DL-ROM neural network parameters associated to two larger  FOM dimensions, $N_h = 128 \times 128 = 16384$ and $N_h = 256 \times 256 = 65536$, by fixing $N = 64$ as rPOD dimension, for a prescribed degree of accuracy $\epsilon_{rel} = 4.03 \times 10^{-3}$. The use of pretraining is possible in this framework because the POD-DL-ROM neural network does not depend on $N_h$, but only on $N$. Pretraining then allows to  reduce the GPU training computational times of a factor $3-7$, as shown in \tablename~\ref{tb:tb2}, where we compare the training times, in presence of pretraining, with the ones of the network trained from scratch.  Irrespectively of pre-training, testing computational times increase with respect to the case $N_h=4096$ due to the dependence on $N_h$ of the matrix-vector product $\mathbf{V}_N\mathbf{u}_N$ required to recover the final POD-DL-ROM approximation, since $\mathbf{V}_N \in \mathbb{R}^{N_h \times N}$ has a higher number of rows.
\begin{table}[h!]
\centering
\begin{tabular}{|l|l|l|l|l|}
\hline
 & $\#$params & $\#$epochs & total time & test $[s]$\\ \hline
POD-DL-ROM ($N_h = 16384$) & 269057 & 1378 & 38 m & 0.06 \\ \hline
POD-DL-ROM PRETRAINED $(N_h = 16384)$ & 269057 & 165 & 5 m & 0.06 \\ \hline
POD-DL-ROM ($N_h = 65536$) & 269057 & 1540 &  42 m & 3 \\ \hline
POD-DL-ROM PRETRAINED ($N_h = 65536$) & 269057 & 461 & 12 m &  3 \\ \hline
\end{tabular}
\caption{\textit{Test 2}: GPU computational times of pretrained and from scratch POD-DL-ROM for $N_h = 16384, 65536$.}
\label{tb:tb2}
\end{table}

Then, we report the results referred to a larger parameter space, in the case $N_h = 4096$ and $N = 64$. In particular, we use the optimal weights associated to $\mathcal{P} = 12.9 \cdot [0.06, 0.2]  \times 12.9 \cdot [0.03,0.1 ]$ $\textnormal{cm\textsuperscript{2}/ms}$ as initial guess of the POD-DL-ROM neural network parameters in the case $\mathcal{P} = 12.9 \cdot [0.02, 0.2]  \times 12.9 \cdot [0.01,0.1 ]$ $\textnormal{cm\textsuperscript{2}/ms}$. We show the GPU training computational times in \tablename~\ref{tb:tb3} for a prescribed level of accuracy, i.e. $\epsilon_{rel} = 4.03 \times 10^{-3}$, and a fixed number of training-parameter instances $N_{train} = 25$, using pretraining and not, respectively. Once again, the use of pretraining allows us to speed up the construction of a POD-DL-ROM remarkably.
\begin{table}[h!]
\centering
\begin{tabular}{|l|l|l|l|l|}
\hline
 & $\#$params & $\#$epochs & total time & test $[s]$\\ \hline
POD-DL-ROM & 269057 & 1486 & 41 m & 0.015 \\ \hline
POD-DL-ROM PRETRAINED & 269057 & 588 & 16 m & 0.015 \\ \hline
\end{tabular}
\caption{\textit{Test 2}: GPU computational times of pretrained and from scratch POD-DL-ROM for $\mathcal{P} = 12.9 \cdot [0.02, 0.2]  \times 12.9 \cdot [0.01,0.1 ]$ $\textnormal{cm\textsuperscript{2}/ms}$.}
\label{tb:tb3}
\end{table}

\subsection{Test 3: nonlinear elastodynamics for hyperelastic compressible materials}

We now consider the solution of an elastodynamics problem, consisting of the following initial/boundary-value problem \cite{gurtin1982anintroduction} for nonlinear elasticity equations, in a three-dimensional beam $\Omega = (0, 1) \times (0, 5) \times (0, 1)$ cm:
\begin{equation}
\left\{
\begin{aligned}
&\rho \frac{\partial^2 \mathbf{d}}{\partial t^2} - \textnormal{div}\left( \mathbf{P}(\mathbf{d}) \right) = \mathbf{f} &\qquad & (\mathbf{x}, t) \in \Omega \times (0, T), \\
&\mathbf{d} = \mathbf{0}& \qquad &(\mathbf{x}, t) \in \Gamma_D \times(0, T), \\
&\mathbf{P}(\mathbf{d}) \mathbf{n} = \mathbf{0} & \qquad & (\mathbf{x}, t) \in \Gamma_N \times ( 0, T), \\
&\mathbf{d}(0) = \mathbf{0} & \qquad & \mathbf{x} \in \Omega,  \ t = 0 \\
& \frac{\partial \mathbf{d}}{\partial t}(0) = \mathbf{0} & \qquad & \mathbf{x} \in \Omega, \ t = 0.
\end{aligned}
\right.
\label{eq:elastodynamics}
\end{equation}
Here we consider $\rho = 1$ kg/cm\textsuperscript{3}, $\mathbf{f} = (-0.01, 0, -0.02)$ kg/(cm s\textsuperscript{2}), $\Gamma_D = \{(x,z) \in (0,1)^2, y = 0 \}$ and $\Gamma_N = \partial \Omega \char`\\ \Gamma_D$.
We consider a St. Venant-Kirchhoff constitutive law involving a hyperelastic nonlinear model to describe the behavior of compressible materials \cite{ogden1997nonlinear}, characterized by the following strain energy function
\begin{equation*}
\psi(\mathbf{F}) = \nu \mathbf{E} : \mathbf{E} + \frac{\lambda}{2} (\text{tr}(\mathbf{E}))^2.
\end{equation*}
Here $\mathbf{E} = \tfrac{1}{2}(\mathbf{F}^T\mathbf{F} - \mathbf{I}) = \tfrac{1}{2} (\mathbf{C} - \mathbf{I})$ is the Green-Lagrange strain tensor, $\mathbf{F} = {\bf I} + \nabla {\bf d}$ is the deformation tensor defined in terms of the displacement ${\bf d} = (d_1, d_2, d_3)$, $\nu$ and $\lambda$ are the Lam\'e coefficients, and 
\begin{equation*}
\mathbf{P}(\mathbf{F}) = \mathbf{F} \left( 2\nu \mathbf{E} + \lambda \text{tr}(\mathbf{E}) \mathbf{I} \right)
\end{equation*}
 is  the first Piola-Kirchhoff stress tensor. Here $n_\mu = 2$ parameters are considered, given by   the Young modulus $\mu_1$ and the Poisson ratio $\mu_2$, belonging to the parameter space $\mathcal{P}= [1,3]\textnormal{ Pa} \times [0.25, 0.42]$;  they affect the expression of the Lam\'e coefficients as follows
\begin{equation*}
\nu = \frac{\mu_1}{2(1+\mu_1)} \quad \textnormal{and} \quad \lambda = \frac{\mu_1\mu_2}{(1+\mu_2)(1-2\mu_2)}.
\end{equation*}
Equations (\ref{eq:elastodynamics}) are discretized in space by means of quadratic $(\mathbb{P}_2)$ finite elements, yielding a dynamical system of dimension $N_h = 5674 \times 3 =  17022$. For time integration, we use the generalized-$\alpha$ method \cite{chung1993atime} over the interval $(0, T)$, with $T = 15$ s and a time-step $\Delta t = 0.2$ s.

We consider $N_t = 75$ time instances over $(0, T)$, $N_{train} = 10 \times 5 = 50$ training-parameter instances $\boldsymbol{\mu}_{train} = (1 + i2/9, 0.25 + j0.0425)$, for $i = 0, \ldots, 9$ and $j = 0, \ldots, 4$, and $N_{test} = 9 \times 4 = 36$ testing-parameter instances $\boldsymbol{\mu}_{test} = (1.111 + i2/9, 0.2712 + j0.0425)$, for $i = 0, \ldots, 9$ and $j = 0, \ldots, 4$. We set the rPOD dimension to \textcolor{black}{$N = 64$ for each of the three components of the displacement} %\times 3 = 192$} 
and the dimension of the nonlinear trial manifold $\tilde{S}_n$ to \textcolor{black}{$n = 3$, for a total number of degrees of freedom of the POD-DL-ROM solution equal to 3.}  The maximum number of epochs is  $N_{epochs} = 10000$, the batch size is $N_b = 20$ and, regarding the early-stopping criterion, we stop the training if the loss function does not decrease \textcolor{black}{along} 500 epochs. The training and testing phases are performed on a GTX 1070 8GB GPU. 

\begin{figure}[b!]
\centering
\includegraphics[scale=0.27]{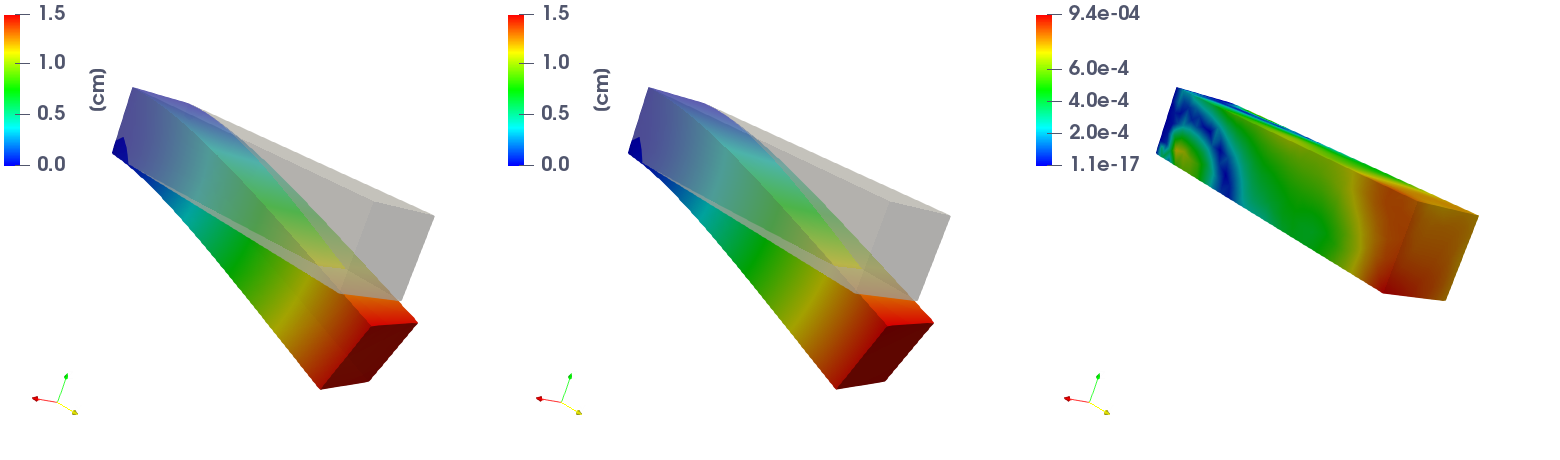}
\caption[\textit{Test 3}: FOM, POD-DL-ROM solutions and relative error $\boldsymbol{\epsilon}_k$, for the testing-parameter instance $\boldsymbol{	\mu}_{test} =(2.88 \textnormal{ Pa}, 0.3987)$ at $T = 15$ s, with $n=3$.]{\textit{Test 3}: FOM (left), POD-DL-ROM (center) solutions and relative error $\boldsymbol{\epsilon}_k$ (right), for the testing-parameter instance $\boldsymbol{	\mu}_{test} =(2.88 \textnormal{ Pa}, 0.3987)$ at $T = 15$ s, with $n=3$.}
\label{fig:Fig5-7}
\end{figure}

\begin{figure}[ht!]
\centering
\includegraphics[scale=0.34]{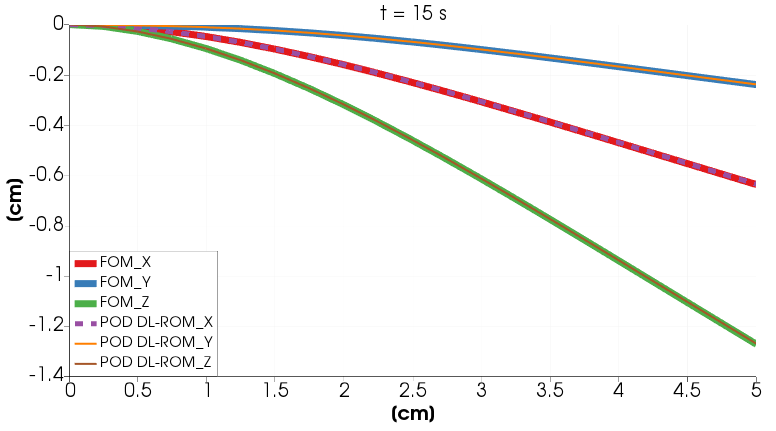}
\caption{\textit{Test 3}: FOM and POD-DL-ROM solutions components over the longitudinal axis, for the testing-parameter instance $\boldsymbol{	\mu}_{test} =(2.88 \textnormal{ Pa}, 0.3987)$ at $T = 15$ s, with $n=3$.}
\label{fig:Fig5-8}
\end{figure}

In \figurename~\ref{fig:Fig5-8} we report the FOM and the POD-DL-ROM solutions, obtained by \textcolor{black}{choosing $n =3$},  and the three components of the  displacement vector ${\bf d} = (d_1, d_2, d_3)$ over the longitudinal axis, i.e. the line which connects the two points $P_1 =(0.5, 0, 0.5)$ cm and $P_2 = (0.5, 5, 0.5)$ cm, for the testing-parameter instance $\boldsymbol{\mu}_{test} =(2.88 \textnormal{ Pa}, 0.3987)$ at $T = 15$ s. We remark that all the three components are accurately captured by the POD-DL-ROM. 
In \figurename~\ref{fig:Fig5-7} we show the FOM solution and the POD-DL-ROM one, with $n = 3$, along with the relative error (\ref{eq:relative_error}), for the testing-parameter instance $\boldsymbol{\mu}_{test} =(2.88 \textnormal{ Pa}, 0.3987)$ at $T = 15$ s. The maximum relative error, which is associated to the portion of the domain undergoing the maximum displacement, is about $10^{-3}$.

Here we also want to investigate how the use of pretraining involving different fidelity models impacts on the POD-DL-ROM technique, \textcolor{black}{starting from the previous low-fidelity model}. In particular, we consider the nearly-incompressible Neo-Hookean constitutive law \cite{ogden1997nonlinear}, an hyperelastic nonlinear model whose strain energy function is specified in terms of an isochoric-volumetric splitting
\begin{equation*}
\psi(\mathbf{F}) = \underbrace{\frac{G}{2}\left( \bar{I}_{1} - 3 \right)}_{\psi_{\text{iso}}} + \underbrace{\frac{K}{4}\left( (J-1)^2 + (\ln J)^2\right)}_{\psi_{\text{vol}}},
\end{equation*}
where $\bar{I}_{1} = J^{-2/3} I_1 = J^{-2/3} \text{tr}(\mathbf{C})$, $J = \text{det}(\mathbf{F})$, $G$ is the shear modulus and $K$ the bulk modulus. The coefficients $G$ and $K$ depend on the Young modulus and the Poisson coefficient and are defined as follow
\begin{equation*}
G = \frac{\mu_1}{2(1 + \lambda)} \quad \textnormal{and} \quad K = \frac{2}{3} G + \lambda.
\end{equation*}
For this model, we consider  $n_\mu = 2$ parameters, belonging to the parameter space $\mathcal{P}= [0.1,1]\textnormal{ Pa} \times [0.3, 0.45]$, and an external force $\mathbf{f} = (-0.001t, 0, -0.002t)$; moreover, the final time is $T =22.5$ s, that is, we enlarge the time interval, and the time-step is set equal to $0.25$ s. We consider $N_t = 90$ time instances over $(0, T)$, $N_{train} = 50$ training-parameter instances and $N_{test} = 36$ testing-parameter instances uniformly distributed over the parameter space. We use the optimal weights and biases found on the first low-fidelity model, as initial guess for the parameters of the POD-DL-ROM on this second configuration \textcolor{black}{which thus features {\em (i)} a more involved constitutive law, {\em (ii)} a different parameter space which reflects in larger displacements, and {\em (iii)} a longer time interval where to compute the dynamics.}

\begin{figure}[ht!]
\vspace{-0.25cm}
\centering
\includegraphics[scale=0.24]{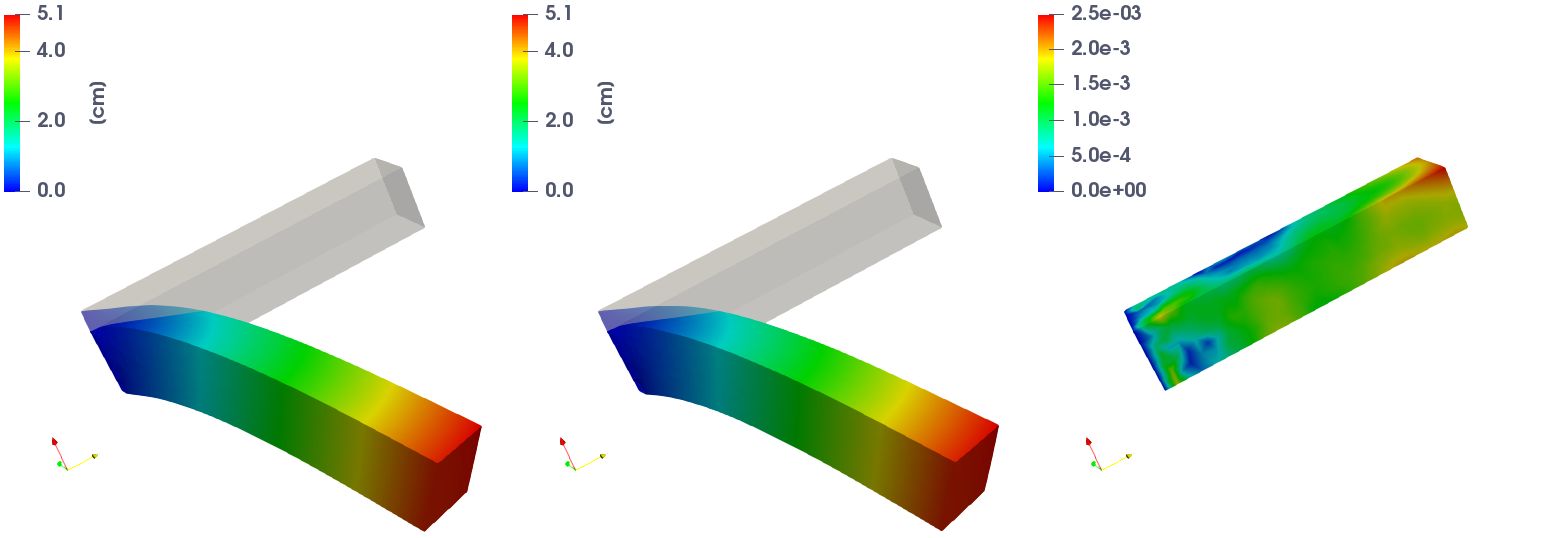}
\includegraphics[scale=0.24]{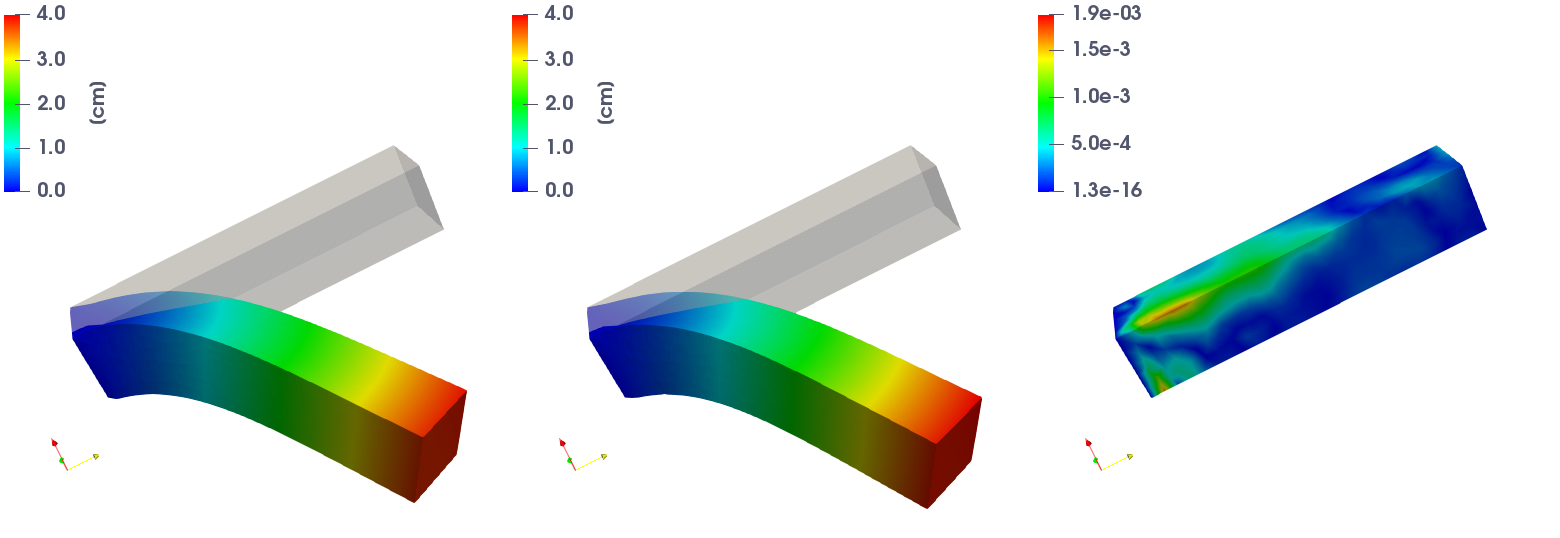}
\caption[\textit{Test 3}: FOM, POD-DL-ROM solutions and relative error $\boldsymbol{\epsilon}_k$, for the testing-parameter instances $\boldsymbol{	\mu}_{test} =(0.25 \textnormal{ Pa}, 0.32)$ and $\boldsymbol{\mu}_{test} =(0.95 \textnormal{ Pa}, 0.43)$ at $T = 22.5$ s, with $n=3$.]{\textit{Test 3}: FOM (left), POD-DL-ROM (center) solutions and relative error $\boldsymbol{\epsilon}_k$ (right), for the testing-parameter instances $\boldsymbol{	\mu}_{test} =(0.25 \textnormal{ Pa}, 0.32)$ (top) and $\boldsymbol{\mu}_{test} =(0.95 \textnormal{ Pa}, 0.43)$ (bottom) at $T = 22.5$ s, with $n=3$.}
\label{fig:Fig5-9}
\end{figure}

\begin{figure}[ht!]
\centering
\includegraphics[scale=0.325]{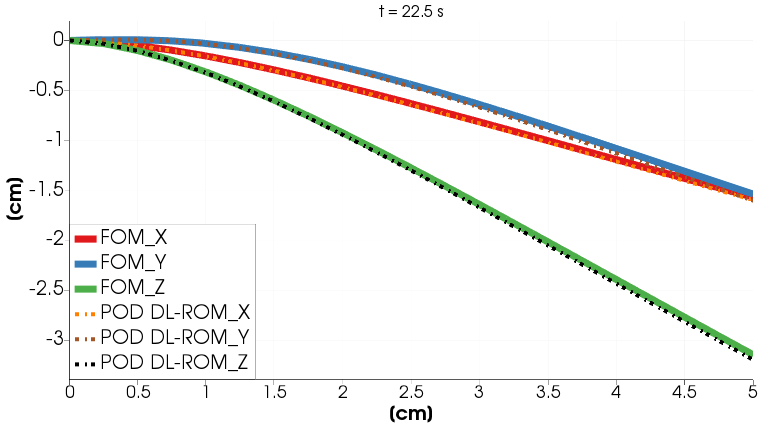}
\caption{\textit{Test 3}: FOM and POD-DL-ROM solutions components over the longitudinal axis, for the testing-parameter instance $\boldsymbol{	\mu}_{test} =(0.25 \textnormal{ Pa}, 0.32)$ at $T = 22.5$ s, with $n=3$.}
\label{fig:Fig5-10}
\vspace{-0.15cm}
\end{figure}

In \figurename~\ref{fig:Fig5-9} we show the FOM and DL-ROM solutions, with $n =3$, together with the relative error (\ref{eq:relative_error}), for the testing-parameter instances $\boldsymbol{	\mu}_{test} =(0.25 \textnormal{ Pa}, 0.32)$ and $\boldsymbol{	\mu}_{test} =(0.95 \textnormal{ Pa}, 0.43)$ at $T = 22.5$ s. In \figurename~\ref{fig:Fig5-10} we compare the FOM and POD-DL-ROM, displaying the three components of the displacement  vector over the longitudinal axis, for the testing-parameter instance $\boldsymbol{\mu}_{test} =(0.25 \textnormal{ Pa}, 0.32)$ at $T = 22.5$ s.

In \tablename~\ref{tb:tb5-4} we finally compare the GPU total and testing computational times of the POD-DL-ROM neural network with and without the use of pretraining. In particular, the use of pretraining allows to strongly reduce the total training and validation time.  The testing computational time, which refers to the time needed by the POD-DL-ROM to compute $N_t = 90$ time instances for a testing-parameter instance, is equal to 0.006 s, and is remarkably lower than the final time $T = 22.5$ s, that is our technique is able to return even faster than real-time solutions. \textcolor{black}{In particular, the testing time reflects in a speed-up, if we consider the time required by the solution of the FOM\footnote{\textcolor{black}{The simulation is performed on} 20 cores of 1.7 TB node (192 Intel\textsuperscript{\textregistered} Xeon Platinum\textsuperscript{\textregistered} 8160 2.1GHz cores) of the HPC cluster available at MOX, Politecnico di Milano.}, equal to $4.12 \times 10^4$.}
\begin{table}[h!]
\centering
\begin{tabular}{|l|l|l|l|l|}
\hline
 & $\#$params & $\#$epochs & total time & test $[s]$\\ \hline
POD-DL-ROM & 270259 & 6490 & 63 m & 0.006 \\ \hline
POD-DL-ROM PRETRAINED & 270259 & 1519 & 15 m & 0.006 \\ \hline
\end{tabular}
\caption{\textit{Test 3}: GPU computational times of pretrained and from scratch POD-DL-ROM.}
\label{tb:tb5-4}
\end{table}

\subsection{Test 4: unsteady Navier-Stokes equations}

\begin{figure}[b!]
\centering
\vspace{-0.15cm}
\includegraphics[scale=0.315]{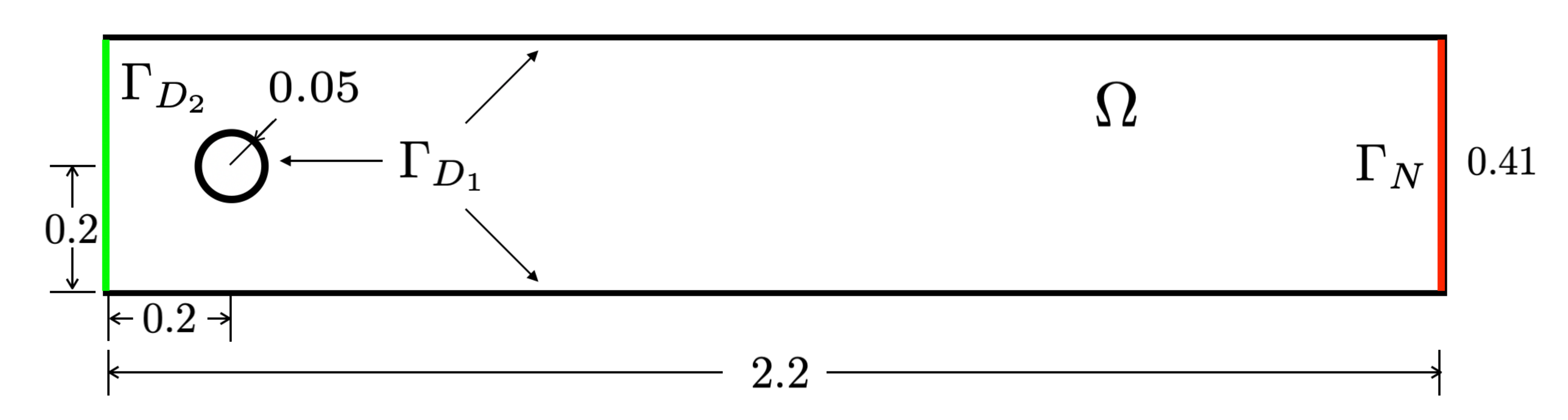}
\caption{\textit{Test 4}: geometrical configuration, domain and boundaries.}
\label{fig:Fig5-11}
\end{figure}

We finally focus on the unsteady Navier-Stokes equations  \cite{quarteroni1994numerical} for incompressible flows in primitive variables (fluid velocity $\mathbf{u}$ and pressure $p$), considering  the flow around a cylinder test case, a well-known benchmark for the evaluation of numerical algorithms for incompressible Navier-Stokes equations in the laminar case:  \\
\begin{equation}
\left\{
\begin{aligned}
& \rho \frac{\partial \mathbf{u}}{\partial t} + \rho \mathbf{u} \cdot \nabla \mathbf{u}   - \nabla \cdot \boldsymbol{\sigma}(\mathbf{u}, p) =  \mathbf{0} & \qquad  & (\mathbf{x},t) \in \Omega \times (0,T), \\
&\nabla \cdot \mathbf{u} = 0  & \qquad & (\mathbf{x},t) \in \Omega\times (0,T),\\
& \mathbf{u} =  \mathbf{0} & \qquad & (\mathbf{x},t) \in \Gamma_{D_1} \times (0,T),\\
& \mathbf{u} =  \mathbf{h} & \qquad & (\mathbf{x},t) \in \Gamma_{D_2} \times (0,T),\\
& \boldsymbol{\sigma}(\mathbf{u}, p) \mathbf{n}= \mathbf{0} & \qquad & (\mathbf{x},t) \in \Gamma_{N}\times (0,T), \\
&  \mathbf{u}(0) = \mathbf{0} & \qquad & \mathbf{x} \in \Omega, \ t = 0.
\end{aligned}
\right.
\label{eq:NS}
\end{equation}
The domain consists in a two-dimensional pipe with a circular obstacle, i.e. $\Omega = (0, 2.2) \times (0, 0.41) \char`\\ \bar{B}_r(0.2,0.2)$ with radius $r=0.05$ (see \figurename~\ref{fig:Fig5-11} for a sketch of the geometry); the boundary is given by $\partial \Omega = \Gamma_{D_1} \cup \Gamma_{D_2} \cup \Gamma_{N}$, where   
$\Gamma_{D_1} = \{ x_1 \in [0, 2.2], x_2 = 0\} \cup \{ x_1 \in [0, 2.2], x_2 = 0.41\} \cup \partial B_{0.05}( (0.2, 0.2) ) $, being $B_r({\bf x}_c)$ the ball of radius $r>0$ centered at ${\bf x}_c$, $\Gamma_{D_2} = \{ x_1 = 0, x_2 \in [0, 0.41] \}$, and  $\Gamma_N = \{x_1 = 2.2, x_2 \in [0, 0.41] \}$, while $\mathbf{n}$ denotes the (outward directed) normal unit vector to $\partial \Omega$. We denote by 
 $\rho$ the fluid density, and by $\boldsymbol{\sigma}$ the stress tensor, 
\begin{equation*}
\boldsymbol{\sigma}(\mathbf{u}, p) = -p \mathbf{I} + 2 \nu \boldsymbol{\epsilon}(\mathbf{u}).
\end{equation*}
Here $\nu$ denotes the dynamic viscosity of the fluid, while $\boldsymbol{\epsilon}(\mathbf{u})$ is the strain tensor, \vspace{-0.1cm}
\begin{equation*}
\boldsymbol{\epsilon}(\mathbf{u}) = \frac{1}{2} \big( \nabla \mathbf{u} + \nabla \mathbf{u} ^T \big). \vspace{-0.1cm}
\end{equation*}
% \cite{NS}).
The density of the fluid is $\rho=1$, no-slip boundary conditions are applied on $\Gamma_1$, a parabolic inflow profile  \begin{equation}
\label{eq:h}
\mathbf{h}(\mathbf{x},t; \mu) = \left( \frac{4 U(t, \mu) x_2 (0.41-x_2)}{0.41^2} , 0 \right), \qquad \mbox{where } \ \  U(t; \mu) = \mu \sin(\pi t / 8),
\end{equation}
is prescribed at the inlet $\Gamma_{D_2}$, while zero-stress Neumann conditions are imposed at the outlet $\Gamma_N$.  We consider as parameter ($n_\mu=1$)  $\mu \in \mathcal{P}=[1,2]$, which reflects on the Reynolds number varying in the range $[66,133]$. Equations (\ref{eq:NS}) have been discretized in space by means of linear-quadratic $(\mathbb{P}_2-\mathbb{P}_1$), inf-sup stable, finite elements, and in time through a BDF of order 2 with semi-implicit treatment of the convective term (see, e.g., \cite{forti2015semi}) over the time interval $T=8$, with a time-step $\Delta t = 2 \times 10^{-3}$.

We uniformly sample $N_t =400$ time instances and consider $N_{train} = 11$ and $N_{test} = 10$ training- and testing-parameter instances uniformly distributed over $\mathcal{P}$. \textcolor{black}{We are interested in reconstructing the velocity field, for which  the FOM dimension is equal to $N_h = 32446 \times 2 = 64892$, selecting $N = 256$ as  dimension of the rPOD basis  for each of the two velocity components. We choose $n = 2$ as  dimension of the nonlinear trial manifold $\tilde{S}_n$.} %for each of the two components of the solution -- once again, we match the dimension $n_{\mu} + 1$ for each component -- for a total number of degrees of freedom of the POD-DL-ROM solution equal to 4.
 We highlight the possibility, by using POD-DL-ROM, to reconstruct the field of interest, i.e the velocity $\mathbf{u}$, without the need of taking into account the approximation of the pressure $p$.

\begin{figure}[ht!]
\vspace{-0.15cm}
\centering
\includegraphics[scale=0.3]{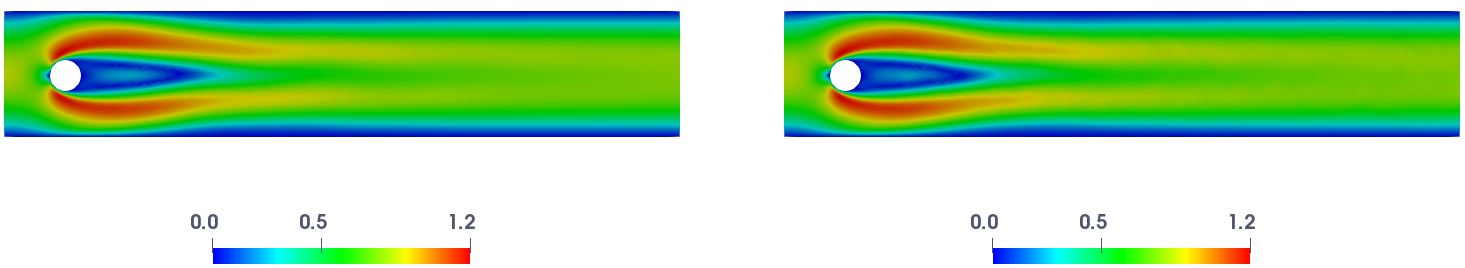} \\
\vspace{0.5cm}
\includegraphics[scale=0.3]{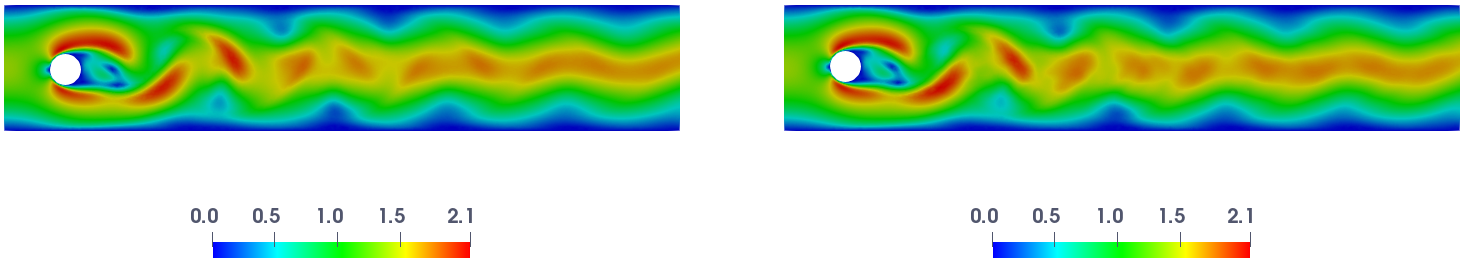} \\
\caption[\textit{Test 4}: FOM and POD-DL-ROM solutions for the testing-parameter instances $\mu_{test} = 1.05$ (Re = 70) and $\mu_{test} = 1.75$ (Re = 117) at $t = 5.64$, with $n=2$.]{\textit{Test 4}: FOM (left) and POD-DL-ROM (right) solutions for the testing-parameter instances $\mu_{test} = 1.05$ (top) and $\mu_{test} = 1.75$ (bottom) at $t = 5.64$, with $n=2$.}
\label{fig:Fig5-12}
\end{figure}

\begin{figure}[ht!]
\centering
\includegraphics[scale=0.35]{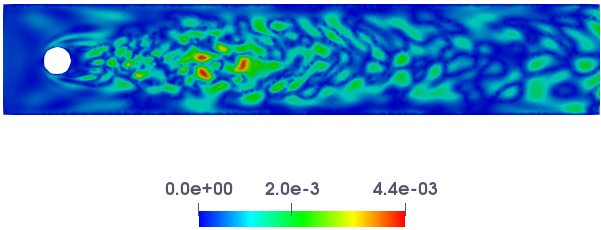} \\
\vspace{0.5cm}
\includegraphics[scale=0.36]{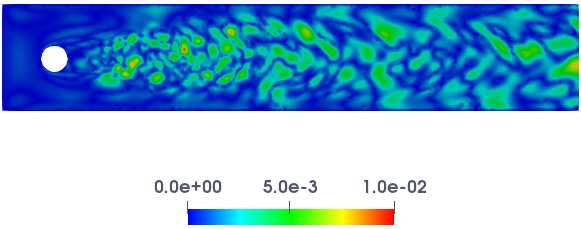}
\caption{\textit{Test 4}: Relative error for the testing-parameter instances $\mu_{test} = 1.05$ (top) and $\mu_{test} = 1.75$ (bottom) at $t = 5.64$.}
\label{fig:Fig5-13}
\vspace{-0.15cm}
\end{figure}

% \times 2 = 512$ 
%
%and the one of the nonlinear trial manifold is $n = 2 \times 2 = 4$. 
In \figurename~\ref{fig:Fig5-12} we compare the FOM and POD-DL-ROM solutions, the latter for $n = 2$, together with the relative error $\boldsymbol{\epsilon}_k$ in \figurename~\ref{fig:Fig5-13}, for two testing-parameter instances $\mu_{test} = 1.05$ (Re = 70) and $\mu_{test} = 1.75$ (Re = 117) at $t = 5.64$. We highlight the  ability of the POD-DL-ROM approximation to accurately capture the variability of the solution over the parameter space $\mathcal{P}$: indeed, in the case Re = 70 (\figurename~\ref{fig:Fig5-12}, top) we do not assist to any vortex shedding; this latter is instead present in the case Re = 117 (\figurename~\ref{fig:Fig5-12}, bottom).

The computational training and testing time of the POD-DL-ROM neural network on a Tesla V100 32 GB GPU are equal to 50 minutes and 0.1 seconds, respectively.

\begin{figure}[ht!]
\vspace{-0.1cm}
\centering
\includegraphics[scale=0.25]{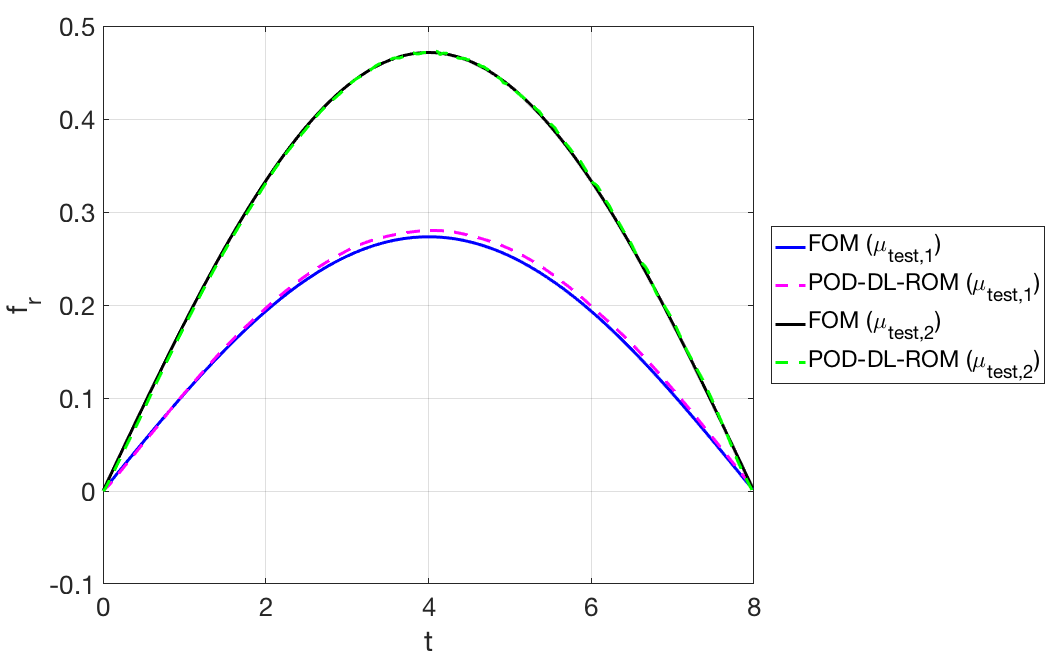} \\
\caption{\textit{Test 4}: FOM and POD-DL-ROM flow rates for the testing-parameter instances $\mu_{test,1} = 1.05$  and $\mu_{test,2} = 1.75$.}
\vspace{-0.15cm}
\label{fig:Fig5-14}
\end{figure}

To show the ability of the POD-DL-ROM to provide accurate evaluations of output quantities of interest, we consider the computation of the flow rate on the outflow boundary $\Gamma_N$, defined as
\begin{equation*}
f_r(t; {\mu}) = \int_{\Gamma_N} {\bf u}(\mathbf{x}, t; {\mu}) \cdot \mathbf{n} \; d\sigma.
\end{equation*}
%where $\mathbf{n}$ is the outward directed normal unit vector to $\Gamma_N$.
 \textcolor{black}{In this respect, in \figurename~\ref{fig:Fig5-14} we show the FOM and POD-DL-ROM flow rates over time, for the two testing-parameter instances $\mu_{test,1} = 1.05$  and $\mu_{test,2} = 1.75$. The POD-DL-ROM technique is able to capture the shape of $f_r$, related to the prescribed $\mu$-dependent input profile in (\ref{eq:h}), in both cases, introducing a maximum relative error equal to 1.65$\%$.}

\textcolor{black}{We then investigate the use of pretraining when aiming at reducing the complexity of the solution of the same problem, whose high-fidelity discretization is set on a finer computational mesh. In this case, we increase   the FOM dimension to $N_h = 128764 \times 2 = 257528$, and train the networks starting from the optimal parameters found on the low-fidelity model related to $N_h = 32446 \times 2 = 64892$}. In particular, 434 epochs, which results in a training computational time of 10 minutes, are required to achieve in the case  $N_h = 128764 \times 2 = 257528$ the same accuracy $\epsilon_{rel} = 1.4 \times 10^{-2}$ obtained for the previous case with  $N_h = 32446 \times 2 = 64892$.
Finally, we show in \tablename~\ref{tb:tb6} the speed-ups introduced, at testing time, by the use of the POD-DL-ROM technique with respect the solution of the FOM\footnote{\textcolor{black}{The simulations are performed on} 20 cores of 1.7 TB node (192 Intel\textsuperscript{\textregistered} Xeon Platinum\textsuperscript{\textregistered} 8160 2.1GHz cores) of the HPC cluster available at MOX, Politecnico di Milano.} \textcolor{black}{when aiming at evaluating the fluid velocity over the interval $(0, T)$, for both cases $N_h = 64892$ and 257528.}
\begin{table}[h!]
\centering
\begin{tabular}{|l|l|l|}
\hline
 & $N_h = 64892$ & $N_h = 257528$ \\ \hline
speed-up & $2.15 \times 10^5$ &  $6.59 \times 10^5$ \\ \hline
\end{tabular}
\caption{\textit{Test 4}: POD-DL-ROM speed-ups.}
\label{tb:tb6}
\end{table}

We remark that ensuring ROM stability in the classical POD-Galerkin framework usually requires additional computational efforts, such as a suitable enrichment of the velocity reduced basis, and a consequent increase of the size of the ROM; see, e.g., \cite{DALSANTO2019186, doi:10.1002/nme.4772, ROZZA20071244}.

\section{Conclusions}

In this work we proposed a strategy to enhance DL-ROMs in order to make the offline training stage dramatically faster. Indeed, a key aspect in the setting of DL-ROMs concerns computational efficiency during the offline (or training) stage, which is also related with the curse of dimensionality. This strategy, which we refer to as POD-DL-ROM, overcomes the main computational bottleneck of the DL-ROM technique, namely the (strong) limitation related to the FOM dimension $N_h$. In particular, it exploits \textit{(i)} dimensionality reduction of FOM snapshots by means of randomized POD (or randomized SVD) and \textit{(ii)} a suitable multi-fidelity pretraining stage exploiting snapshots computed through lower-fidelity models to initialize the parameters of neural networks in a sequential procedure. Moreover, the POD-DL-ROM approximations retain all the features of DL-ROM solutions, enabling extremely efficient testing computational times. 
% Indeed, we have assessed the computational performance of POD-DL-ROMs in {\em (i)} a linear advection-diffusion-reaction problem,  {\em (ii)} a nonlinear diffusion-reaction problem arising from cardiac electrophysiology, {\em (iii)} structural mechanics, and
%{\em (iv)} fluid dynamics. As a result, POD-DL-ROMs are shown to yield extremely efficient numerical approximations to (scalar and vector) nonlinear time-dependent parametrized PDEs, thus providing a {\em turn-key} strategy to build ROMs only relying on a set of FOM snapshots, and ultimately leading to the possibility to simulate in more than real-time, during the {\em online} testing stage, physical phenomena occurring on a time scale of seconds.

We assessed computational performance, numerical accuracy and robustness of the POD-DL-ROM technique on several time-dependent parametrized PDEs, namely {\em (i)} a linear advection-diffusion-reaction problem,  {\em (ii)} a nonlinear diffusion-reaction problem arising from cardiac electrophysiology, {\em (iii)} nonlinear elastodynamics for hyperelastic compressible materials, and {\em (iv)} fluid dynamics.  In all these cases, POD-DL-ROMs are able to match the intrinsic dimension of the problems investigated, to overcome the main computational bottleneck shown by conventional projection-based methods, and to make the training phase of ROMs extremely fast. Through the numerical test cases assessed in Section \ref{sec:num_res}, POD-DL-ROMs have shown to yield extremely efficient numerical approximations to (scalar and vector) nonlinear time-dependent parametrized PDEs, thus providing a {\em turn-key} strategy to build ROMs only relying on a set of FOM snapshots, and ultimately leading to the possibility to simulate in more than real-time, during the {\em online} testing stage, physical phenomena occurring on a time scale of seconds.

\section*{Acknowledgments}

%The authors have been partially supported by the ERC Advanced Grant iHEART, ``An integrated heart model for the simulation of the cardiac function'', 2017-2022, P.I. A. Quarteroni (ERC2016AdG, project ID: 740132). Moreover, 
We gratefully acknowledge Prof. A. Quarteroni and Prof. L. Dede' (MOX, Politecnico di Milano) for their insightful discussions and  useful remarks.

%, ultimately leading to the possibility to solve in more than real-time, during the {online} testing stage, parametrized PDEs modeling physical phenomena whose time scale is  seconds.} \\

%\section*{References}

\bibliography{references}

\end{document}